\documentclass[11pt,a4paper,twoside,final]{scrartcl}
\usepackage{a4wide}
\usepackage{amsfonts}
\usepackage{amsmath}
\usepackage{amssymb}
\usepackage{amsthm} 
\usepackage{algorithm}
\usepackage{algorithmic}
\usepackage{pgfplots}
\pgfplotsset{compat=newest}
\usepackage{graphicx}
\usepackage{color}
\usepackage{colortbl}
\usepackage{showkeys}
\usepackage{enumerate}
\usepackage{subcaption}
\usepackage{verbatim}
\usepackage{bm}
\usepackage[bookmarks=false]{hyperref}
\usepackage{tikz}
\usepackage{tikz-cd}
\pgfplotsset{colormap={whitered}{color(0cm)=(white);
	color(1cm)=(orange!75!red)}
}

\newtheorem{theorem}{Theorem}[section]
\newtheorem{lemma}[theorem]{Lemma}
\newtheorem{remark}[theorem]{Remark}
\newtheorem{generalisation}[theorem]{Generalisation}
\newtheorem{definition}[theorem]{Definition}
\newtheorem{example}[theorem]{Example}
\newtheorem{corollary}[theorem]{Corollary}
\newtheorem{proposition}[theorem]{Proposition}

\newenvironment{Example}{\goodbreak \begin{example}\slshape}{\bend\end{example}}

\makeatletter
\def\imod#1{\allowbreak\mkern10mu({\operator@font mod}\,\,#1)}
\makeatother

\makeatletter 
 
\@addtoreset{algorithm}{chapter} 
\makeatother

\numberwithin{equation}{section}
\numberwithin{table}{section}
\numberwithin{figure}{section}

\newcommand{\bend}{\hspace*{0ex} \hfill \hbox{\vrule height
    1.5ex\vbox{\hrule width 1.4ex \vskip 1.4ex\hrule  width 1.4ex}\vrule
    height 1.5ex}}

\long\def\symbolfootnote[#1]#2{\begingroup \def\thefootnote{\fnsymbol{footnote}}\footnote[#1]{#2}\endgroup}
\makeatletter
\let\@fnsymbol\@arabic
\makeatother

\usepackage{tabularx}
\newcolumntype{L}[1]{>{\raggedright\arraybackslash}p{#1}} \newcolumntype{C}[1]{>{\centering\arraybackslash}p{#1}} \newcolumntype{R}[1]{>{\raggedleft\arraybackslash}p{#1}}

\title{A note on transformed Fourier systems for the approximation of non-periodic signals}

\date{}\author{
	Robert Nasdala\thanks{
    	Faculty of Mathematics, Chemnitz University of Technology, D-09107 Chemnitz, Germany.\newline E-mail:                \href{mailto:robert.nasdala@math.tu-chemnitz.de}{robert.nasdala@math.tu-chemnitz.de}
	}
	\and
	Daniel Potts\thanks{
		Faculty of Mathematics, Chemnitz University of Technology, D-09107 Chemnitz, Germany.\newline E-mail:                \href{mailto:potts@math.tu-chemnitz.de}{potts@math.tu-chemnitz.de}
	}
}

\hypersetup{pdfauthor={Robert Nasdala, Daniel Potts},
		pdftitle={Efficient multivariate approximation on the cube},
		pdfsubject={},
		pdfcreator = {pdflatex},
		plainpages=false,
		pdfstartview=FitH,       
		pdfview=FitH,            
		pdfpagemode=UseOutlines, 
		bookmarksnumbered=true, 
		bookmarksopen=false,     
		bookmarksopenlevel=0,   
		colorlinks=true,       
		linkcolor=black,
		citecolor=black,
		urlcolor=black}
\begin{document}

\maketitle

\medskip

\begin{abstract}
	A variety of techniques have been developed for the approximation of non-periodic functions. 
	In particular, there are approximation techniques based on rank-$1$ lattices and transformed {rank-$1$} lattices, including methods that use sampling sets consisting of Chebyshev- and tent-transformed nodes.
	We compare these methods with a parameterized transformed Fourier system that yields similar $\ell_2$-approximation errors.
\end{abstract}

\section{Introduction}
For the approximation of non-periodic functions defined on the cube $[0,1]^d$, fast algorithms based on Chebyshev- and tent-transformed rank-$1$ lattice methods have been introduced and studied in \cite{DiNuPi14,PoVo15,CoKuNuSu16,SuNuCo14,IrrKrPi18,GoSzYo19,KuoMiNoNu19}.
Recently, we suggested a general framework for transformed rank-$1$ lattice approximation, in which functions defined on a cube $[0,1]^d$ (or on $\mathbb{R}^d$) are periodized onto the torus $\mathbb{T}^d\simeq[0,1)^d$, \cite{NaPo18,NaPo19}.
In these approaches we define parameterized families $\psi(\cdot,\bm\eta):\left[0,1\right]^d\to\left[0,1\right]^d, \bm\eta\in\mathbb{R}_{+}^{d}$ of transformations that, depending on the parameter choice, yield a certain smoothening effect when composed with a given non-periodic function.
This periodization strategies also lead to general parameterized classes of orthonormal systems in weighted Hilbert spaces.
However, these methods have the natural drawback of singularities appearing at the boundary points of the cube, so that any approximation error estimates have to be done with respect to weighted $L_{\infty}$- and $L_{2}$-norms.

We summarize some crucial properties of rank-$1$ lattice approximation.
Then, we compare the approximation with a half-periodic cosine system and tent-transformed sampling nodes \cite{IsNo08,iserles2006high,AdcockDiss,Adcock11,SuNuCo14},
the Chebyshev approximation \cite{PoVo15,KuoMiNoNu19}, 
as well as the general framework for the parameterized transformed Fourier system \cite{NaPo19}.
We discuss numerical results in up to dimension $d=7$ and highlight the controlled smoothening effect when varying the parameter $\bm\eta$ in the transformed Fourier systems.

\section{Approximation methods}
At first, we summarize the main ideas of the Fourier approximation with sampling sets in the form of rank-$1$ lattices \cite{SLKA87,CoKuNu10,kaemmererdiss}.
Afterwards, we consider Chebyshev- and tent-transformed {rank-$1$} lattices in the context of Chebyshev and cosine approximation methods \cite{SuNuCo14,KuoMiNoNu19}.
Finally, we outline the transformed Fourier system for the approximation of non-periodic signals, as introduced in \cite{NaPo19}, and provide two examples of parameterized transformations.

\subsection{Fourier approximation}
For any frequency set $I \subset \mathbb{Z}^d$ of finite cardinality $|I|<\infty$ we denote the space of all multivariate trigonometric polynomials supported on $I$ by
\begin{align*} 
	\Pi_{I} := \mathrm{span}\left\{ \mathrm{e}^{2\pi\mathrm i \mathbf k \cdot \mathbf x} = \prod_{\ell=1}^{d}\mathrm{e}^{2\pi\mathrm i k_\ell x_\ell} : \mathbf k = (k_{1},\ldots,k_{d})^{\top}\in I, \mathbf x = (x_{1},\ldots,x_{d})^{\top}\in \mathbb{T}^d\right\}.
\end{align*}
Trigonometric polynomials are orthonormal with respect to the $L_2(\mathbb T^d)$-scalar product
\begin{align*}
	(f,g)_{L_2(\mathbb T^d)} := \int_{\mathbb T^d} f(\mathbf x)\, \overline{g(\mathbf x)}\,\mathrm d\mathbf x, 
	\quad f,g\in L_2(\mathbb{T}^d).
\end{align*}
For all $\mathbf k\in \mathbb{Z}^d$ we denote the \textsl{Fourier coefficients} $\hat h_{\mathbf k}$ by
\begin{align*}
	\hat h_{\mathbf k} 
	:= (h, \mathrm{e}^{2\pi\mathrm i \mathbf k (\cdot)})_{L_2(\mathbb T^d)} 
	= \int_{\mathbb T^d} h(\mathbf x)\,\mathrm{e}^{-2\pi\mathrm i \mathbf k \cdot \mathbf x} \,\mathrm d\mathbf x,
\end{align*}
and the corresponding \textsl{Fourier partial sum} by $S_{I}h(\mathbf x) := \sum_{\mathbf k\in I} \hat h_{\mathbf k}\,\mathrm{e}^{2\pi\mathrm i \mathbf k\mathbf \cdot \mathbf x}$.

We use sampling nodes in a \textsl{rank-$1$ lattice} $\Lambda(\mathbf z, M)$ of size $M\in\mathbb{N}$ generated by the vector $\mathbf z\in\mathbb{Z}^d$, that is defined as
\begin{align}\label{def:R1L}
	\Lambda(\mathbf z, M) 
	:= \left\{ \mathbf x_{j}^{\mathrm{latt}} := \frac{j}{M}\mathbf z \mod \mathbf 1 \in\mathbb{T}^{d} : j = 0,\ldots, M-1 \right\},
\end{align}
with $\mathbf 1 := (1,\ldots,1)^{\top}$, which allows the fast evaluation of Fourier partial sums via \cite[Algorithm~3.1]{kaemmererdiss}.
For any frequency set $I \subset \mathbb{Z}^d$ the \textsl{difference set} is given by
\begin{align}\label{def:difference_set}
	\mathcal{D}(I) &:= \{\mathbf{k} \in \mathbb{Z}^d : \mathbf{k}=\mathbf{k}_1-\mathbf{k}_2 \text{ with } \mathbf{k}_1,\mathbf{k}_2 \in I \}.
\end{align}
We define the \textsl{reconstructing {rank-$1$} lattice} $\Lambda(\mathbf{z}, M, I)$ as a {rank-$1$} lattice $\Lambda(\mathbf{z}, M)$ for which the condition
\begin{align}\label{eq:recon_R1L}
	\mathbf{t} \cdot \mathbf{z} \not\equiv 0 \,(\bmod{M}) \quad \text{for all } \mathbf{t}\in\mathcal{D}(I)\setminus\{\mathbf 0\} 
\end{align}
holds.
Given a reconstructing {rank-$1$} lattice $\Lambda(\mathbf z,M,I)$, we have exact integration for all multivariate trigonometric polynomials $p \in \Pi_{\mathcal{D}(I)}$, see \cite{SLKA87}, so that
\begin{align}\label{eq:exact_integration_property}
	\int_{\mathbb{T}^d} p(\mathbf x) \,\mathrm{d}\mathbf x 
	= \frac{1}{M} \sum_{j =0}^{M-1} p(\mathbf x_j), 
	\quad \mathbf x_j\in\Lambda(\mathbf z,M,I).
\end{align}
In particular, for $h \in \Pi_{I}$ and $\mathbf k\in I$ we have $h(\cdot)\,\mathrm{e}^{-2\pi\mathrm{i}\mathbf k(\cdot)} \in \Pi_{\mathcal{D}(I)}$ and
\begin{align} 
	\label{eq:exact_integration_formula}
	\hat h_{\mathbf k}
	= \int_{\mathbb{T}^d} h(\mathbf x) \,\mathrm{e}^{-2\pi\mathrm{i}\mathbf k\cdot\mathbf x} \,\mathrm{d}\mathbf x 
	= \frac{1}{M} \sum_{j =0}^{M-1} h(\mathbf x_j) \,\mathrm{e}^{-2\pi\mathrm{i}\mathbf k\cdot\mathbf x_j}, 
	\quad \mathbf x_j\in\Lambda(\mathbf z,M,I).
\end{align}
Next, we focus on functions in the Wiener algebra $\mathcal{A}(\mathbb{T}^d)$ containing all $L_{1}(\mathbb{T}^d)$-functions with absolutely summable Fourier coefficients $\hat{h}_{\mathbf k}$ given by
\begin{align}\label{def:WienerAlgebra}
	\mathcal{A}(\mathbb{T}^d) := \left\{h \in L_{1}(\mathbb{T}^d) : \sum_{k\in\mathbb{Z}^d} |\hat{h}_{\mathbf k}| < \infty \right\}.
\end{align}
For an arbitrary function ${h \in \mathcal{A}(\mathbb{T}^d)\cap\mathcal{C}(\mathbb{T}^d)}$ and lattice points ${\mathbf x_j \in \Lambda(\mathbf z,M,I)}$ we lose the former mentioned exact integration property and get \textsl{approximated Fourier coefficients} $\hat h_{\mathbf k}^\Lambda$ of the form
\begin{align*} 
	\hat h_{\mathbf k}
	&\approx \hat h_{\mathbf k}^\Lambda
	:= \frac{1}{M} \sum_{j =0}^{M-1} h(\mathbf x_j)\,\mathrm{e}^{-2\pi\mathrm i \mathbf k\cdot\mathbf x_j}
\end{align*}
leading to the \textsl{approximated Fourier partial sum} $S_{I}^{\Lambda} h$ given by
\begin{align*}
	h(\mathbf x)
	&\approx S_{I}^{\Lambda} h(\mathbf x)
	:= \sum_{\mathbf k\in I} \hat h_{\mathbf k}^\Lambda \,\mathrm{e}^{2\pi\mathrm i \mathbf k\cdot\mathbf x}. 
\end{align*}
For the matrix-vector-expression with respect to the frequency set $I_{\mathrm{latt}}\subset\mathbb{Z}^d$ we put 
\begin{align*}
	\mathbf F_{\mathrm{latt}} 
	:= \left\{ 
		\mathrm{e}^{2\pi\mathrm{i}\mathbf k\cdot \mathbf x_{j}^{\mathrm{latt}}}
	\right\}_{j=0, \mathbf k\in I_{\mathrm{latt}} }^{M-1}, 
	\quad
	\mathbf h_{\mathrm{latt}} := \left(h( \mathbf x_{j}^{\mathrm{latt}} )\right)_{j=0}^{M-1}.
\end{align*}
The evaluation of the function $h$ and the reconstruction of the approximated Fourier coefficients $\mathbf{\hat h} := (\hat{h}_{\mathbf k}^{\Lambda})_{\mathbf k \in I_{\mathrm{latt}}}$ are realized by the fast Algorithms outlined in~\cite[Algorithm~3.1 and 3.2]{kaemmererdiss} that solve the systems
\begin{align} \label{eq:Fourier_system_MatVec}
	\mathbf{h_{\mathrm{latt}}}
	=
	\mathbf F_{\mathrm{latt}} \mathbf{\hat h}
	\quad \text{and} \quad 
	\mathbf{\hat h}
	=
	\frac{1}{M}\mathbf F_{\mathrm{latt}}^{*} \mathbf{h_{\mathrm{latt}}},
\end{align}
where we have $\mathbf F_{\mathrm{latt}}^{*}\mathbf F_{\mathrm{latt}} = M \mathbf I$ by construction with the identity matrix $\mathbf I\in\mathbb{C}^{|I_{\mathrm{latt}}|\times |I_{\mathrm{latt}}|}$. 
\subsection{Cosine approximation}
Next, we consider the half-periodic cosine system 
\begin{align}\label{def:cosine_system}
	\left\{ \lambda_{\mathbf k}(\mathbf x) := \sqrt{2}^{\|\mathbf k\|_{0}}\prod_{\ell=1}^{d}\cos(\pi k_\ell x_\ell) \right\}_{\mathbf k\in I_{\mathrm{tent}}},
	I_{\mathrm{tent}} \subset \mathbb{N}_{0}^d, \mathbf x \in [0,1]^d,
\end{align}
with the zero-norm $\|\mathbf k\|_0 := |\{\ell\in\{1,\ldots,d\}:k_\ell\neq 0\}|$ and $\sqrt{2}^{\|\mathbf k\|_0} := \prod_{\ell=1}^{d} \sqrt{2}^{\|k_\ell\|_0}$.
In \cite{IsNo08} it is pointed out that this system can alternatively be defined in one dimension over the domain $t\in [-1,1]$ as the system $\lambda_0(x) = \frac{1}{\sqrt{2}}, \lambda_k(t) = \cos(k\pi t), \tilde\lambda_k(t) = \sin((k-\frac{1}{2})\pi t)$, which yields the original cosine system after applying the transformation $t=2x-1$.

The cosine system~\eqref{def:cosine_system} is orthonormal with respect to the $L_{2}([0,1]^d)$-scalar product given by
\begin{align*}
	(f,g)_{L_2([0,1]^d)} := \int_{[0,1]^d} f(\mathbf x)\, \overline{g(\mathbf x)}\,\mathrm d\mathbf x, 
	\quad f,g\in L_2([0,1]^d).
\end{align*}
For $\mathbf k \in \mathbb{Z}^d$ the cosine coefficient
\index{cosine coefficient!$\hat{h}_{\mathbf k}^{\mathrm{cos}}$}  of a function $h\in L_2\left([0,1]^d\right)$ is naturally defined as
$\hat{h}_{\mathbf k}^{\mathrm{cos}} := \left( h, \lambda_{\mathbf k} \right)_{L_{2}\left([0,1]^d\right)}$ and for $I\subset\mathbb{Z}^d$ the corresponding cosine partial sum is given by 
$S_{I} h(\mathbf x) := \sum_{\mathbf k\in I} \hat{h}_{\mathbf k}^{\mathrm{cos}} \,\lambda_{\mathbf k}(\mathbf x)$.
We transfer the crucial properties of the Fourier system via the tent transformation 
\begin{align} \label{def:tent_trafo}
	\psi(\mathbf x) := (\psi_1(x_1),\ldots,\psi_d(x_d))^{\top}, \quad
	\psi_\ell(x_\ell) = \begin{cases}
		2x_\ell & \text{for } 0 \leq x_\ell < \frac{1}{2}, \\
		2 - 2x_\ell & \text{for } \frac{1}{2} \leq x_\ell \leq 1.
	\end{cases}
\end{align}
We have sampling nodes in the tent-transformed rank-$1$ lattice $\Lambda_{\psi}(\mathbf z, M)$ defined as
\begin{align}\label{def:R1L_tenttransformed}
	\Lambda_{\psi}(\mathbf z, M)
	:= \left\{ \mathbf y_{j}^{\mathrm{cos}} := \psi\left(\mathbf x_{j}^{\mathrm{latt}}\right) : \mathbf x_{j}^{\mathrm{latt}}\in\Lambda(\mathbf z,M), j = 0,\ldots, M-1 \right\}
\end{align}
and we speak of a reconstructing tent-transformed rank-$1$ lattice $\Lambda_{\psi}(\mathbf z,M,I)$ if the underlying rank-$1$ lattice is a reconstructing one.
Recalling the definition of difference sets $\mathcal{D}(I)$ in \eqref{def:difference_set},
multivariate trigonometric polynomials $h(\cdot), h(\cdot)$ and $\lambda_{\mathbf k}(\cdot)$ that are in $\Pi_{\mathcal{D}(I)}$ and supported on $\mathbf k\in I\subset\mathbb{N}_0^d$ inherit the exact integration property~\eqref{eq:exact_integration_property}, because with the tent transformation as in \eqref{def:tent_trafo} and transformed nodes $\mathbf y_j^{\mathrm{cos}} = \psi(\mathbf x_j^{\mathrm{latt}})\in\Lambda_{\psi}(\mathbf z,M,I)$ with $\mathbf x_j^{\mathrm{latt}} = (x_{1}^{j},\ldots,x_{d}^{j})^{\top}\in\Lambda(\mathbf z,M,I)$ we have
\begin{align*} 
	\hat{h}_{\mathbf k}^{\mathrm{cos}}
	= \int_{[0,1]^d} h(\mathbf y) \,\lambda_{\mathbf k}(\mathbf y) \,\mathrm{d}\mathbf y 
	&= \sqrt{2}^{\|\mathbf k\|_{0}}\int_{\mathbb{T}^d} h(\psi(\mathbf x)) \,\prod_{\ell=1}^{d}\cos(2\pi k_\ell x_\ell) \,\mathrm{d}\mathbf x \\
	&= \frac{\sqrt{2}^{\|\mathbf k\|_{0}}}{2^d}\int_{\mathbb{T}^d} h(\psi(\mathbf x)) \,\left(\mathrm{e}^{2\pi\mathrm{i}\mathbf k\cdot \mathbf x} + \mathrm{e}^{-2\pi\mathrm{i}\mathbf k\cdot \mathbf x}\right) \,\mathrm{d}\mathbf x \\
	&= \frac{\sqrt{2}^{\|\mathbf k\|_{0}}}{2^d} \frac{1}{M}\sum_{j=0}^{M-1} h(\psi(\mathbf x_j)) \,\left(\mathrm{e}^{2\pi\mathrm{i}\mathbf k\cdot \mathbf x_j} + \mathrm{e}^{-2\pi\mathrm{i}\mathbf k\cdot \mathbf x_j}\right) \\
	&= \sqrt{2}^{\|\mathbf k\|_{0}} \frac{1}{M}\sum_{j=0}^{M-1} h(\psi(\mathbf x_j)) \,\prod_{\ell=1}^{d}\cos(2\pi k_\ell x_{\ell}^{j}) \\
	&= \frac{1}{M} \sum_{j =0}^{M-1} h(\mathbf y_j^{\mathrm{cos}}) \,\lambda_{\mathbf k}(\mathbf y_j^{\mathrm{cos}}).
\end{align*}
For an arbitrary function ${h \in \mathcal{C}\left([0,1]^d\right)}$, we lose the former mentioned exactness and define the \textsl{approximated cosine coefficients}
$\hat h_{\mathbf k}^{\mathrm{cos},\Lambda}$ of the form
\begin{align*} 
	\hat{h}_{\mathbf k}^{\mathrm{cos}}
	\approx \hat{h}_{\mathbf k}^{\mathrm{cos},\Lambda}
	:= \frac{1}{M} \sum_{j =0}^{M-1} h(\mathbf y_j^{\mathrm{cos}})\,\lambda_{\mathbf k}(\mathbf y_j^{\mathrm{cos}}),
	\quad\mathbf y_j^{\mathrm{cos}} \in \Lambda_{\psi}(\mathbf z,M,I),
\end{align*}
and obtain \textsl{approximated cosine partial sum} $S_{I}^{\Lambda} h$ given by
\begin{align}\label{def:approx_cosine_Partsum}
	h(\mathbf x)
	&\approx S_{I}^{\Lambda} h(\mathbf x)
	:= \sum_{\mathbf k\in I} \hat{h}_{\mathbf k}^{\mathrm{cos},\Lambda} \,\lambda_{\mathbf k}(\mathbf x). 
\end{align}
In matrix-vector-notation we have
\begin{align*}
	\mathbf C
	:= \left\{
		\lambda_{\mathbf k} \left( \mathbf y_{j}^{\mathrm{cos}} \right) 
	\right\}_{j=0, \mathbf k\in I_{\mathrm{tent}}}^{M-1},
	\quad 
	\mathbf h_{\mathrm{tent}} 
	:= \left(h( \mathbf y_{j}^{\mathrm{cos}} )\right)_{j=0}^{M-1}.
\end{align*}
Both the evaluation of $h$ and the reconstruction of the approximated cosine coefficients ${\mathbf{\hat{h}} := \left\{ \hat{h}_{\mathbf k}^{\mathrm{cos},\Lambda} \right\}_{\mathbf k\in I_{\mathrm{tent}}}}$ is realized by solving the systems
\begin{align}\label{eq:Tent_system_MatVec}
	\mathbf{h_{\mathrm{tent}}}
	=
	\mathbf C \mathbf{\hat h}
	\quad\text{and}\quad
	\mathbf{\hat h}
	=
	\frac{1}{M}\mathbf C^{*} \mathbf{h_{\mathrm{tent}}},
\end{align}
where we have $\mathbf C^{*}\mathbf C = M \mathbf I$ by construction with the identity matrix $\mathbf I\in\mathbb{C}^{|I_{\mathrm{tent}}|\times |I_{\mathrm{tent}}|}$.
Fast algorithms for solving both systems are described in \cite{SuNuCo14,KuoMiNoNu19}.

\subsection{Chebyshev approximation}
We consider the Chebyshev system, that is defined for $\mathbf x \in[0,1]^d$ and a finite frequency set $I_{\mathrm{cheb}}\subset \mathbb{N}_{0}^d$ as
\begin{align}\label{def:T_k}
	\left\{ 
		T_{\mathbf k}(\mathbf x) 
		:= \sqrt{2}^{\|\mathbf k\|_0} \prod_{\ell=1}^{d}  \cos\left(k_{\ell} \arccos(2x_{\ell}-1)\right) 
	\right\}_{\mathbf k\in I_{\mathrm{cheb}}}.
\end{align}
The Chebyshev system \eqref{def:T_k} is an orthonormal system with respect to the weighted scalar product
\begin{align*}
	(T_{\mathbf k_1}, T_{\mathbf k_2})_{L_{2}([0,1]^d,\omega)} 
	:= \int_{[0,1]^d} T_{\mathbf k_1}(\mathbf x) \, T_{\mathbf k_2}(\mathbf x) \, \omega(\mathbf x) \,\mathrm{d}\mathbf x,
	\quad
	\omega(\mathbf x) := \prod_{\ell=1}^{d}\frac{2}{\pi \sqrt{4x_\ell(1-x_\ell)}}.
\end{align*}
The \textsl{Chebyshev coefficients} of a function $h\in L_{2}\left([0,1]^d,\omega\right)$ are naturally defined as
$\hat{h}_{\mathbf k}^{\mathrm{cheb}} := (h, T_{\mathbf k})_{L_{2}\left([0,1]^d,\omega\right)}, \mathbf k\in\mathbb{Z}^d$ and for $I\subset\mathbb{Z}^d$ the corresponding Chebyshev partial sum is given by 
$S_{I}h(\mathbf x) := \sum_{\mathbf k\in I} \hat{h}_{\mathbf k}^{\mathrm{cheb}} \,T_{\mathbf k}(\mathbf x)$.
\index{Chebyshev partial sum!$S_{I}h$}
We transfer some properties of the Fourier system via the \textsl{Chebyshev transformation}
\begin{align}\label{def:Cheb_trafo}
	\psi(\mathbf x) := (\psi_1(x_1),\ldots,\psi_d(x_d))^{\top},
	\quad
	\psi_\ell(x_\ell) := \frac{1}{2} + \frac{1}{2}\cos\left( 2\pi \left(x_\ell-\frac{1}{2}\right) \right), \quad x_\ell\in\left[0,1\right].
\end{align}
We have sampling nodes in the Chebyshev-transformed rank-$1$ lattice $\Lambda_{\psi}(\mathbf z, M)$ defined as
\begin{align}\label{def:R1L_chebtransformed}
	\Lambda_{\psi}(\mathbf z, M)
	:= \left\{ 
		\mathbf y_{j}^{\mathrm{cheb}} := \psi\left(\mathbf x_{j}^{\mathrm{latt}}\right) : \mathbf x_{j}^{\mathrm{latt}}\in\Lambda(\mathbf z,M), j = 0,\ldots, M-1 
	\right\}.
\end{align}
It inherits the reconstruction property~\eqref{eq:recon_R1L} of the underlying reconstructing rank-$1$ lattice $\Lambda(\mathbf z,M,I)$ and is denoted by $\Lambda_{\psi}(\mathbf z, M, I)$.
We note that Chebyshev transformed sampling nodes are fundamentally connected to Padua points and Lissajous curves, as well as certain interpolation methods that are outlined in \cite{BoCaMaViXu06, DeEr15}. 

Recalling the definition of difference sets $\mathcal{D}(I)$ in \eqref{def:difference_set}, 
multivariate trigonometric polynomials $h(\cdot)$ and $h(\cdot)\,T_{\mathbf k}(\cdot)$ are in $\Pi_{\mathcal{D}(I)}$ and supported on $\mathbf k\in I\subset\mathbb{N}_0^d$ inherit the exact integration property~\eqref{eq:exact_integration_property}, because with the Chebyshev transformation $\psi$ as in \eqref{def:Cheb_trafo} and transformed nodes $\mathbf y_j^{\mathrm{cheb}} = \psi(\mathbf x_j^{\mathrm{latt}})\in\Lambda_{\psi}(\mathbf z,M,I)$ with $\mathbf x_j^{\mathrm{latt}} = (x_{1}^{j},\ldots,x_{d}^{j})^{\top}\in\Lambda(\mathbf z,M,I)$ we have
\begin{align*} 
	\hat{h}_{\mathbf k}^{\mathrm{cheb}}
	= \int_{[0,1]^d} h(\mathbf y) \,T_{\mathbf k}(\mathbf y)\,\omega(\mathbf y) \,\mathrm{d}\mathbf y 
	&= \sqrt{2}^{\|\mathbf k\|_{0}}\int_{\mathbb{T}^d} h(\psi(\mathbf x)) \,\prod_{\ell=1}^{d}\cos(2\pi k_\ell x_\ell) \,\mathrm{d}\mathbf x \\
	&= \frac{\sqrt{2}^{\|\mathbf k\|_{0}}}{2^d}\int_{\mathbb{T}^d} h(\psi(\mathbf x)) \,\left(\mathrm{e}^{2\pi\mathrm{i}\mathbf k\cdot \mathbf x} + \mathrm{e}^{-2\pi\mathrm{i}\mathbf k\cdot \mathbf x}\right) \,\mathrm{d}\mathbf x \\
	&= \frac{\sqrt{2}^{\|\mathbf k\|_{0}}}{2^d} \frac{1}{M}\sum_{j=0}^{M-1} h(\psi(\mathbf x_j)) \,\left(\mathrm{e}^{2\pi\mathrm{i}\mathbf k\cdot \mathbf x_j} + \mathrm{e}^{-2\pi\mathrm{i}\mathbf k\cdot \mathbf x_j}\right) \\
	&= \sqrt{2}^{\|\mathbf k\|_{0}} \frac{1}{M}\sum_{j=0}^{M-1} h(\psi(\mathbf x_j)) \,\prod_{\ell=1}^{d}\cos(2\pi k_\ell x_{\ell}^{j}) \\
	&= \frac{1}{M} \sum_{j =0}^{M-1} h(\mathbf y_j^{\mathrm{cheb}}) \,T_{\mathbf k}(\mathbf y_j^{\mathrm{cheb}}).
\end{align*}
For an arbitrary function $h \in L\left([0,1]^d,\omega\right)\cap\mathcal{C}\left([0,1]^d\right)$, we lose the former mentioned exactness and define the \textsl{approximated Chebyshev coefficients} 
\index{Chebyshev coefficient!approximated, $\hat{h}_{\mathbf k}^{\mathrm{cheb},\Lambda}$}
$\hat{h}_{\mathbf k}^{\mathrm{cheb},\Lambda}$ of the form
\begin{align*} 
	\hat{h}_{\mathbf k}^{\mathrm{cheb}}
	\approx \hat{h}_{\mathbf k}^{\mathrm{cheb},\Lambda}
	:= \frac{1}{M} \sum_{j =0}^{M-1} h(\mathbf y_j^{\mathrm{cheb}})\,T_{\mathbf k}(\mathbf y_j^{\mathrm{cheb}}),
	\quad\mathbf y_j^{\mathrm{cheb}} \in \Lambda_{\psi}(\mathbf z,M,I),
\end{align*}
leading to the \textsl{approximated Chebyshev partial sum}
\begin{align}\label{def:approx_Cheb_Partsum}
	h(\mathbf x)
	&\approx S_{I}^{\Lambda} h(\mathbf x)
	:= \sum_{\mathbf k\in I} \hat{h}_{\mathbf k}^{\mathrm{cheb},\Lambda}  \,T_{\mathbf k}(\mathbf x). 
\end{align}
In matrix-vector-notation this reads as
\begin{align*}
	\mathbf T
	:= \left\{
		T_{\mathbf k}(\mathbf y_{j}^{\mathrm{cheb}})
	\right\}_{j=0, \mathbf k\in I_{\mathrm{cheb}}}^{M-1},
	\quad
	\mathbf h_{\mathrm{cheb}} := \left(h( \mathbf y_{j}^{\mathrm{cheb}} )\right)_{j=0}^{M-1}.
\end{align*}
The evaluation of $h$ as well as the reconstruction of the approximated Chebyshev coefficients ${ \mathbf{\hat h} := \left(\hat{h}_{\mathbf k}^{\mathrm{cheb},\Lambda}\right)_{\mathbf k \in I_{\mathrm{cheb}}} }$ of $h$ are realized by fast Algorithms outlined in~\cite{PoVo15,SuNuCo14, KuoMiNoNu19}, that solve the systems
\begin{align}\label{eq:Chebyshev_system_MatVec}
	\mathbf{h_{\mathrm{cheb}}}
	=
	\mathbf T \mathbf{\hat h}
	\quad\text{and}\quad
	\mathbf{\hat h}
	=
	\frac{1}{M}\mathbf T^{*} \mathbf{h_{\mathrm{cheb}}},
\end{align}
where we have $\mathbf T^{*}\mathbf T = M \mathbf I$ by construction with the identity matrix $\mathbf I\in\mathbb{C}^{|I_{\mathrm{cheb}}|\times |I_{\mathrm{cheb}}|}$.

\subsection{Transformed Fourier approximation}
We recall the ideas of a particular family of parameterized torus-to-cube transformations as suggested in \cite{NaPo19}, that generalize the construction idea of the Chebyshev system in composing a mapping with a multiple of its inverse. 

We call a continuously differentiable, strictly increasing mapping $\tilde{\psi}:(0,1) \to \mathbb{R}$ with $\tilde{\psi}(x+\frac{1}{2})$ being odd and $\tilde{\psi}(x) \to \pm\infty$ for $x \to \{0,1\}$ a \textsl{torus-to-$\mathbb{R}$ transformation}.
We obtain a parameterized \textsl{torus-to-cube transformation} ${\psi(\cdot,\eta):[0,1] \to [0,1]}$ with $\eta\in\mathbb{R}_{+} := (0,\infty)$ by putting
\begin{align} \label{def:Trafo_comb}
	\psi(x, \eta) := 
	\begin{cases}
		0 & \text{for}\quad x=0, \\
		\tilde{\psi}^{-1}(\eta \,\tilde{\psi}(x)) & \text{for}\quad x\in\left(0,1\right), \\
		1 & \text{for}\quad x=1,
	\end{cases}
\end{align}
which are continuously differentiable, increasing and have a first derivative $\psi'(\cdot,\eta)\in \mathcal{C}_{0}([0,1])$,
where $\mathcal{C}_{0}\left([0,1]\right)$ denotes the space of all continuous functions vanishing to $0$ towards their boundary points.
It holds
${\psi^{-1}(y,\eta) = \psi\left(y,\frac{1}{\eta}\right)}$ and 
we call ${\varrho(y,\eta) := (\psi^{-1})'(y,\eta) = \psi'\left(y,\frac{1}{\eta}\right)}$ the \textsl{density of $\psi$}. 
In multiple dimensions $d\in\mathbb{N}$ with $\bm \eta = (\eta_1,\ldots,\eta_d)^{\top}$ we put
\begin{align} \label{def:Trafo_comb_mult}
	\psi(\mathbf x, \bm \eta) &:= ( \psi_1(x_1, \eta_1),\ldots,\psi_d(x_d, \eta_d) )^{\top},\\
	\psi^{-1}(\mathbf y, \bm\eta) &:= ( \psi_1^{-1}(y_1,\eta_1),\ldots,\psi_d^{-1}(y_d,\eta_d) )^{\top}, \nonumber\\
	\varrho(\mathbf y,\bm\eta) &:= \prod_{\ell=1}^d \varrho_\ell(y_\ell,\eta_\ell)
	\quad \text{with} \quad 
	\varrho_\ell(y_\ell,\eta_\ell) := \frac{1}{\psi'(\psi^{-1}(y_\ell,\eta_\ell))}, \nonumber
\end{align}
where the univariate torus-to-cube transformations $\psi_\ell(\cdot, \eta_\ell)$ and their corresponding densities $\varrho_\ell(\cdot, \eta_\ell)$ may be different in each coordinate $\ell\in\{1,\ldots,d\}$.

We consider integrable weight functions 
\begin{align*}
	\omega(\mathbf y) := \prod_{\ell=1}^{d}\omega_\ell(y_\ell), \quad \mathbf y\in [0,1]^d,
\end{align*}
such that for any given torus-to-cube transformation $\psi(\cdot,\bm\eta)$ as in \eqref{def:Trafo_comb_mult} we have
\begin{align*}
	\omega(\psi_\ell(\cdot,\eta_\ell)) \psi'(\cdot,\eta_\ell) \in \mathcal{C}_{0}\left([0,1]\right).
\end{align*}
Applying a torus-to-cube transformation to a function $h\in L_{2}([0,1]^d,\omega)\cap\mathcal{C}([0,1]^d)$ generates a periodic function $f\in L_{2}(\mathbb{T}^d)$ of the form
\begin{align}\label{eq:periodization_idea}
	f(\mathbf x) := h(\psi(\mathbf x,\bm\eta))\sqrt{\omega(\psi(\mathbf x,\bm\eta)) \prod_{\ell=1}^{d} \psi'_\ell(x_\ell)}
	\quad \text{with} \quad 
	\|h\|_{L_{2}([0,1]^d,\omega)} = \|f\|_{L_{2}(\mathbb{T}^d)},
\end{align}
that is approximated by the classical Fourier system.
To construct an approximant for the original function $h$ we apply the inverse torus-to-cube transformation to the Fourier system, yielding for a fixed $\bm\eta\in\mathbb{R}_+^d$ the \textsl{transformed Fourier system}
\begin{align} \label{def:trafo_Fou_system}
	\left\{ \varphi_{\mathbf k}(\cdot) := \sqrt{\frac{\varrho(\cdot,\bm\eta)}{\omega(\cdot)}}\,\mathrm{e}^{2\pi\mathrm{i}\mathbf k\cdot\psi^{-1}(\cdot,\bm\eta)} \right\}_{\mathbf k\in I},
\end{align}
which forms an orthonormal system with respect to the weighted $L_{2}\left([0,1]^d,\omega\right)$-scalar product.
For all $\mathbf k\in \mathbb{Z}^d$ the \textsl{transformed Fourier coefficients} $\hat h_{\mathbf k}$ are naturally defined as
\begin{align*}
	\hat h_{\mathbf k} 
	:= (h, \varphi_{\mathbf k})_{L_{2}([0,1]^d,\omega)} 
	= \int_{\left[0,1\right]^d} h(\mathbf y)\,\overline{\varphi_{\mathbf k}(\mathbf y)} \, \omega(\mathbf y) \,\mathrm d\mathbf y,
\end{align*}
and the corresponding \textsl{Fourier partial sum} is given by $S_{I}h(\mathbf y) := \sum_{\mathbf k\in I} \hat h_{\mathbf k}\,\varphi(\mathbf y)$.
The corresponding sampling nodes will be taken from the torus-to-cube-transformed (abbreviated: ttc) rank-$1$ lattice $\Lambda_{\psi}(\mathbf z, M)$ defined as
\begin{align}
	\Lambda_{\psi}(\mathbf z, M)
	:= \left\{ \mathbf y_{j}^{\mathrm{ttc}} := \psi\left(\mathbf x_{j}^{\mathrm{latt}},\bm\eta\right) : \mathbf x_{j}^{\mathrm{latt}}\in\Lambda(\mathbf z,M), j = 0,\ldots, M-1 \right\}
\end{align}
and we speak of a reconstructing torus-to-cube-transformed rank-$1$ lattice $\Lambda_{\psi}(\mathbf z,M,I)$ if the underlying rank-$1$ lattice is a reconstructing one.

Furthermore, the multivariate transformed trigonometric polynomials supported on $I\subset\mathbb{Z}^d$ are given by $\Pi_{I}^{\mathrm{ttc}} := \mathrm{span}\{ \varphi_{\mathbf k} : \mathbf k\in I\}$ and inherit the exact integration property~\eqref{eq:exact_integration_formula}, thus, for $h\in\Pi_{I}^{\mathrm{ttc}}$ we have
\begin{align*} 
\hat h_{\mathbf k}
	= \int_{[0,1]^d} h(\mathbf x) \,\varphi_{\mathbf k}(\mathbf x) \,\mathrm{d}\mathbf x 
	= \frac{1}{M} \sum_{j =0}^{M-1} h(\mathbf y_j^{\mathrm{ttc}}) \,\varphi_{\mathbf k}(\mathbf y_{j}^{\mathrm{ttc}}), 
	\quad \mathbf y_j^{\mathrm{ttc}}\in\Lambda_{\psi}(\mathbf z,M,I).
\end{align*}
For an arbitrary function $h\in L_2([0,1]^d,\omega)\cap\mathcal{C}([0,1]^d)$ we lose the former mentioned exactness and define approximated transformed coefficients of the form
\begin{align*}
	\hat{h}_{\mathbf k}^{\Lambda} 
	:= \frac{1}{M} \sum_{j =0}^{M-1} h(\mathbf y_j^{\mathrm{ttc}}) \,\varphi_{\mathbf k}(\mathbf y_{j}^{\mathrm{ttc}})
\end{align*}
and leads to the \textsl{approximated transformed Fourier partial sum} $S_{I}^{\Lambda} h$ given by
\begin{align}\label{def:approx_trafo_Fou_Partsum}
	h(\mathbf y)
	&\approx S_{I}^{\Lambda} h(\mathbf y)
	:= \sum_{\mathbf k\in I} \hat h_{\mathbf k}^\Lambda \,\varphi_{\mathbf k}(\mathbf y). 
\end{align}
In matrix-vector-notation we have
\begin{align*}
	\mathbf h_{\mathrm{ttc}} 
	:= \left(h( \mathbf y_{j}^{\mathrm{ttc}} )\right)_{j=0}^{M-1},
	\quad
	\mathbf F_{\mathrm{ttc}} 
	:= \left\{
		\varphi_{\mathbf k} \left( \mathbf y_{j}^{\mathrm{ttc}}\right) 
	\right\}_{j=0, \mathbf k\in I_{\mathrm{ttc}}}^{M-1}.
\end{align*}
The evaluation of $h$ and the reconstruction of the approximated transformed Fourier coefficients $\mathbf{\hat{h}} := \left\{ \hat{h}_{\mathbf k}^{\Lambda} \right\}_{\mathbf k\in I_{\mathrm{ttc}}}$ is realized by solving the systems
\begin{align}\label{eq:Trafo_Exp_system_MatVec}
	\mathbf{h_{\mathrm{ttc}}}
	=
	\mathbf F_{\mathrm{ttc}} \mathbf{\hat h}.
	\quad\text{and}\quad
	\mathbf{\hat h}
	=
	\frac{1}{M}\mathbf F_{\mathrm{ttc}}^{*} \mathbf{h_{\mathrm{ttc}}}.
\end{align}
Fast algorithms for solving both systems are described in \cite{NaPo19}.

\subsection{Comparison of the orthonormal systems}
The previously presented approximation approaches are based on very different orthonormal systems and use differently transformed sampling sets, which is summarized in dimension $d=1$ in Table~\ref{table:CompONS} with the definition of the hyperbolic cross $I_{N}^{1}$ given in \eqref{def:HC}.

Given an univariate continuous function $h\in\mathcal{C}([0,1])$, both composition with the tent transformation \eqref{def:tent_trafo} and the Chebyshev transformation \eqref{def:Cheb_trafo} can be interpreted as mirroring a compressed version of $h$ at the point $\frac{1}{2}$, so that $h(\psi(x)) = h(\psi(1-x))$ for all $x\in[0,\frac{1}{2}]$.
In contrast to the the Chebyshev transformation case, for the tent transformation we generally won't expect the resulting function $h\circ\psi$ to be smooth at the point $\frac{1}{2}$, which will be reflected in the approximation results later on.

The parametrized torus-to-cube transformations~\eqref{def:Trafo_comb} are a fundamentally different transformation class in the sense that the periodization effect is caused primarily by the multiplication of $h(\psi(\cdot,\eta))$ with the first derivative $\psi(\cdot,\eta)\in\mathcal{C}_{0}([0,1])$ (assuming a constant weight function $\omega\equiv 1$), so that the function $h(\psi(\cdot,\eta))\sqrt{\psi(\cdot,\eta)}$ ends up being continuously extendable to the torus $\mathbb{T}$.
Additionally, now there is the parameter $\eta$ involved which controls the smoothening effect on the periodized function, see \cite{NaPo19}.

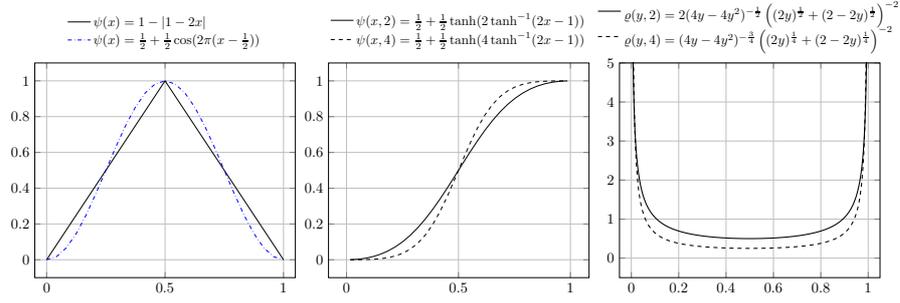
\begin{figure}[t]
	\centering
	\begin{minipage}[b]{.24\linewidth}
	\centering
		\begin{tikzpicture}[
			scale=0.5,
			declare function = {
				tent(\x) 
				= (2*\x)*(0 <= \x)*(\x < 1/2) 
				+ (2 - 2*\x)*(1/2 <= \x)*(\x <= 1) ;
			}
		]
		\begin{axis}[
grid=major,
			xmin=-0.05,	xmax=1.05,ymin=-0.1,ymax=1.1,
			xtick={0,0.5,1.0},
legend style={at={(0.5,1.05)}, anchor=south,legend columns=1,legend cell align=left, font=\small, legend style={draw=none}},
		]
			\addplot[domain=0:1] {tent(x)};
			\addlegendentry{$\psi(x) = 1 - |1-2x|$};
			\addplot[dashdotted, blue, domain=0:1] {0.5+0.5*cos(2*pi*deg(x-1/2))}; 
			\addlegendentry{$\psi(x) = \frac{1}{2} + \frac{1}{2}\cos(2\pi (x-\frac{1}{2}))$};
\end{axis}
	\end{tikzpicture}
	\end{minipage}
\begin{minipage}[b]{.24\linewidth}
	\centering
		\begin{tikzpicture}[
			scale=0.5, declare function = {
				arctanh(\x) = 0.5*(ln(1+\x)-ln(1-\x));
			}]
			\begin{axis}[
				title = {$\psi(x) = $},
				grid = major,
				xtick={0,0.5,1.0},
				legend style={at={(0.5,1.05)}, anchor=south,legend columns=1,legend cell align=left, font=\small, legend style={draw=none}},
			]
\addplot[samples=300] {0.5 + 0.5*tanh(2*arctanh(2*x-1))}; 
			\addlegendentry{$\psi(x,2) = \frac{1}{2} + \frac{1}{2}\tanh(2 \tanh^{-1}(2x-1))$};
			\addplot[dashed, samples=300] {0.5 + 0.5*tanh(3*arctanh(2*x-1))}; 
			\addlegendentry{$\psi(x,4) = \frac{1}{2} + \frac{1}{2}\tanh(4 \tanh^{-1}(2x-1))$};
			\end{axis}
		\end{tikzpicture}
	\end{minipage}
	\begin{minipage}[b]{.24\linewidth}
	\centering
		\begin{tikzpicture}[scale=0.5]
				\pgfmathsetmacro\M{2}
				\pgfmathsetmacro\MM{4}
				\begin{axis}[
					samples=500,  
					xmin=-0.05, xmax=1.05, ymin=-0.5, ymax=5,
grid = major,
					every axis plot/.append style={thick},
					legend style={at={(0.5,1.05)}, anchor=south,legend columns=1,legend cell align=left, font=\small, legend style={draw=none}}
				]
\addplot[black, domain=-0.001:0.999] { 4*(1/\M)*((4*x-4*x^2)^((1/\M)-1))/(((2*x)^(1/\M)+(2-2*x)^(1/\M))^2) };
		\addlegendentry{$\varrho(y,2) = 2 (4y-4y^2)^{-\frac{1}{2}}\left((2y)^\frac{1}{2}+(2-2y)^\frac{1}{2}\right)^{-2}$};
		\addplot[black, dashed, domain=-0.001:0.999] { 4*(1/\MM)*((4*x-4*x^2)^((1/\MM)-1))/(((2*x)^(1/\MM)+(2-2*x)^(1/\MM))^2) };
		\addlegendentry{$\varrho(y,4) = (4y-4y^2)^{-\frac{3}{4}} \left((2y)^\frac{1}{4}+(2-2y)^\frac{1}{4}\right)^{-2}$};
		\end{axis}
		\end{tikzpicture}
	\end{minipage}
	\caption{Left: The tent-transformation \eqref{def:tent_trafo} the Chebyshev-transformation \eqref{def:Cheb_trafo}. 
		Center and right: The parameterized logarithmic transformation \eqref{eq:logarithmic_trafo} and its density function for $\eta\in\{2,4\}$. 
	}
	\label{fig:Trafo_comparison}
\end{figure}

\begin{Example}
	We find various suggestions for torus-to-$\mathbb{R}$ transformations in \cite[Section~17.6]{boyd00}, \cite[Section~7.5]{ShTaWa11} and \cite{NaPo18}.
	We list some induced combined transformations $\psi(x,\eta)$ and the corresponding density function $\varrho(y,\eta) = (\psi^{-1})'(y,\eta)$ in the sense of definition~\eqref{def:Trafo_comb}:
	\begin{itemize}
	\item 
		the logarithmic torus-to-cube transformation
		\begin{align}\label{eq:logarithmic_trafo}
			\psi(x,\eta) := \frac{1}{2} + \frac{1}{2} \tanh (\eta \tanh^{-1}(2x-1)),
			\quad \varrho(y,\eta) =  \frac{4}{\eta}\frac{(4y-4y^2)^{\frac{1}{\eta}-1}}{\left((2y)^\frac{1}{\eta}+(2-2y)^\frac{1}{\eta}\right)^2},
		\end{align}
		based on the mapping
		\begin{align*}
			\tilde{\psi}(x) 
			&= \frac{1}{2}\log\left(\frac{2x}{2-2x}\right) 
			= \tanh^{-1}(2x-1),
		\end{align*}
	\item	
		the error function torus-to-cube transformation
		\begin{align} \label{eq:erf_trafo}
			\psi(x,\eta) = \frac{1}{2}\,\mathrm{erf}(\eta\,\mathrm{erf}^{-1}(2x-1))+\frac{1}{2} ,
			\quad \varrho(y,\eta) = \frac{1}{\eta}\,\mathrm{e}^{(1-\frac{1}{\eta^2})(\mathrm{erf}^{-1}(2y-1))^2},
		\end{align}
		based on the mapping
		\begin{align*}
			\tilde{\psi}(x) &= \mathrm{erf}^{-1}(2x-1),
		\end{align*}
		which is the inverse of the error function 
		\begin{align*}
			\mathrm{erf}(y) = \frac{1}{\sqrt{\pi}} \int_{-y}^{y} \mathrm{e}^{-t^2} \,\mathrm{d}t, \quad y\in\mathbb{R},
		\end{align*}
	\end{itemize}
\end{Example}
In Figure~\ref{fig:Trafo_comparison} we provide a side-by-side comparison of all the previously mentioned transformation mappings.

\begin{table}
\begin{tabular}{|l|C{2.0cm}|C{4cm}|c|}
	\hline
	orthonormal system $\{\varphi_k(x)\}_{k\in I}$ & scalar product weight $\omega$ & sampling transformation $\psi$ & frequency set $I$ \\
	\hline
	$\sqrt{2}^{\|k\|_{0}}\cos(\pi k x)$ & $1$ & $\begin{cases} 2x & \text{for } 0 \leq x < \frac{1}{2}, \\ 2 - 2x & \text{for } \frac{1}{2} \leq x \leq 1.	\end{cases}$ & $I_{N}^{d}\cap\mathbb{N}_{0}^{d}$ \\
	$\sqrt{2}^{\|k\|_{0}} \cos\left(k \arccos(2x-1)\right) $ & $\frac{1}{2\pi \sqrt{x(1-x)}}$ & $\frac{1}{2}+\frac{1}{2}\cos\left( 2\pi (x-\frac{1}{2}) \right)$ & $I_{N}^{d}\cap\mathbb{N}_{0}^{d}$ \\
	$\sqrt{\frac{\varrho(x,\eta)}{\omega(x)}}\,\mathrm{e}^{2\pi\mathrm{i}k\psi^{-1}(x,\eta)}$ & $\omega(x)$ & $\psi(x,\eta)$ & $I_{N}^{d}$ \\
	\hline
\end{tabular}
\caption{Comparison of the univariate orthonormal system, sampling sets and frequency sets from the Chebyshev, cosine and transformed Fourier approximation methods.}
\label{table:CompONS}
\end{table}

\section{Approximation results and error analysis}
Based on the weight function 
\begin{align*}
	\omega_{\mathrm{hc}}(\mathbf k) 
	:= \prod_{\ell=1}^{d} \max(1,|k_j|), 
	\quad \mathbf k\in\mathbb{Z}^d,
\end{align*}
we define the hyperbolic cross index set
\begin{align}\label{def:HC}
	I_{N}^{d} := \left\{ \mathbf k\in\mathbb{Z}^d  : \omega_{\mathrm{hc}}(\mathbf k) \leq N \right\}
\end{align}
and for $\beta\ge0$ we furthermore have the Hilbert spaces
\begin{align}
	\label{def:HbetaRaum}
	\mathcal{H}^{\beta}(\mathbb T^d) 
	:= \left\{ f \in L_2(\mathbb{T}^d) : \|f\|_{\mathcal{H}^{\beta}(\mathbb{T}^d)} := \left( \sum_{\mathbf k\in\mathbb{Z}^d} \omega_{\mathrm{hc}}(\mathbf k)^{2\beta} |\hat f_{\mathbf k}|^2 \right)^{\frac{1}{2}} < \infty \right\}
\end{align}
that are closely related to the Wiener Algebra $\mathcal{A}(\mathbb T^d)$ given in \eqref{def:WienerAlgebra}.
For $\lambda>\frac{1}{2}$ and fixed $d\in\mathbb{N}$ the continuous embeddings $\mathcal{H}^{\beta+\lambda}(\mathbb{T}^d) \hookrightarrow \mathcal{A}(\mathbb T^d)$ was shown in \cite[Lemma~2.2]{KaPoVo13}.	
Next, we introduce the analogue on the cube $[0,1]^d$ for the Hilbert space $\mathcal{H}^{\beta}(\mathbb T^d)$ as in \eqref{def:HbetaRaum}.
We define the space of weighted $L_{2}\left([0,1]^d,\omega\right)$-functions with square summable Fourier coefficients $\hat h_{\mathbf k} := (h, \varphi_{\mathbf k})_{L_2\left([0,1]^d, \omega \right)}$ by
\begin{align}\label{def:Sobolev_space_cube}
\mathcal{H}^{\beta}\left([0,1]^d,\omega\right) 
	:= \left\{ h \in L_2\left([0,1]^d, \omega\right) : 
	\|h\|_{\mathcal{H}^{\beta}\left([0,1]^d,\omega\right)} := \left( \sum_{\mathbf k\in\mathbb{Z}^d} \omega_{\mathrm{hc}}(\mathbf k)^{2\beta} |\hat h_{\mathbf k}|^2 \right)^{\frac{1}{2}} < \infty \right\}.
\end{align}
In case of a constant weight function $\omega \equiv 1$ we just write $\mathcal{H}^{\beta}\left([0,1]^d\right)$.

\begin{figure}[t]
	\centering
\begin{minipage}[c]{0.3\textwidth}
		\centering
		\begin{tikzpicture}[baseline,scale=0.5]
		\begin{axis}[
		grid=major,
		xmin=-0.55,	xmax=0.55,ymin=-0.1,ymax=1.0,
		title = {$h(y) = B_2(y)$},
		xtick={-0.5,0,0.5},
		xticklabels={0,0.5,1},
]
\addplot[smooth] coordinates {
			(-0.5000, 0.7500) (-0.4900, 0.7499) (-0.4800, 0.7496) (-0.4700, 0.7491) (-0.4600, 0.7484) (-0.4500, 0.7475) (-0.4400, 0.7464) (-0.4300, 0.7451) (-0.4200, 0.7436) (-0.4100, 0.7419) (-0.4000, 0.7400) (-0.3900, 0.7379) (-0.3800, 0.7356) (-0.3700, 0.7331) (-0.3600, 0.7304) (-0.3500, 0.7275) (-0.3400, 0.7244) (-0.3300, 0.7211) (-0.3200, 0.7176) (-0.3100, 0.7139) (-0.3000, 0.7100) (-0.2900, 0.7059) (-0.2800, 0.7016) (-0.2700, 0.6971) (-0.2600, 0.6924) (-0.2500, 0.6875) (-0.2400, 0.6824) (-0.2300, 0.6771) (-0.2200, 0.6716) (-0.2100, 0.6659) (-0.2000, 0.6600) (-0.1900, 0.6539) (-0.1800, 0.6476) (-0.1700, 0.6411) (-0.1600, 0.6344) (-0.1500, 0.6275) (-0.1400, 0.6204) (-0.1300, 0.6131) (-0.1200, 0.6056) (-0.1100, 0.5979) (-0.1000, 0.5900) (-0.0900, 0.5819) (-0.0800, 0.5736) (-0.0700, 0.5651) (-0.0600, 0.5564) (-0.0500, 0.5475) (-0.0400, 0.5384) (-0.0300, 0.5291) (-0.0200, 0.5196) (-0.0100, 0.5099) (0.0000, 0.5000) (0.0100, 0.4900) (0.0200, 0.4802) (0.0300, 0.4705) (0.0400, 0.4608) (0.0500, 0.4513) (0.0600, 0.4418) (0.0700, 0.4325) (0.0800, 0.4232) (0.0900, 0.4141) (0.1000, 0.4050) (0.1100, 0.3961) (0.1200, 0.3872) (0.1300, 0.3785) (0.1400, 0.3698) (0.1500, 0.3613) (0.1600, 0.3528) (0.1700, 0.3445) (0.1800, 0.3362) (0.1900, 0.3281) (0.2000, 0.3200) (0.2100, 0.3120) (0.2200, 0.3042) (0.2300, 0.2964) (0.2400, 0.2888) (0.2500, 0.2812) (0.2600, 0.2738) (0.2700, 0.2664) (0.2800, 0.2592) (0.2900, 0.2520) (0.3000, 0.2450) (0.3100, 0.2380) (0.3200, 0.2312) (0.3300, 0.2245) (0.3400, 0.2178) (0.3500, 0.2112) (0.3600, 0.2048) (0.3700, 0.1985) (0.3800, 0.1922) (0.3900, 0.1860) (0.4000, 0.1800) (0.4100, 0.1740) (0.4200, 0.1682) (0.4300, 0.1624) (0.4400, 0.1568) (0.4500, 0.1512) (0.4600, 0.1458) (0.4700, 0.1405) (0.4800, 0.1352) (0.4900, 0.1300) (0.5000, 0.1250)
		};
		\end{axis}
		\end{tikzpicture}
	\end{minipage}
\begin{minipage}[c]{.3\linewidth}
	\centering
		\begin{tikzpicture}[
			scale=0.6,
			declare function = {
				B2univ(\x) 
				= (-\x^2+(3/4))*(0 <= \x)*(\x < 1/2) 
				+ ((1/2)*(\x^2 -3*\x +9/4))*(1/2 <= \x)*(\x <= 1) ;
			},
			declare function = {
				B2mult(\x,\y) 
				= B2univ(\x) * B2univ(\y);
			}
		]
		\begin{axis}[view={45}{30},
			grid = both,
			title = {$h_2(y_1,y_2) = B_2(y_1)\, B_2(y_2)$},
			samples = 20,
			zmin = 0, zmax = 1,
			legend style={at={(0.5,1.05)}, anchor=south,legend columns=2,legend cell align=left, font=\small}
		]
			\addplot3[surf,domain=0:1, domain y=0:1, opacity=0.5] {B2mult(x,y)};
\end{axis}
		\end{tikzpicture}
	\end{minipage}
	\caption{The univariate B-spline $h(y) = B_2(y)$ and the two-dimensional tensored B-spline $h_1(y_1,y_2) = B_2(y_1)\,B_2(y_2)$.
	}
	\label{fig:B2plot}
\end{figure}
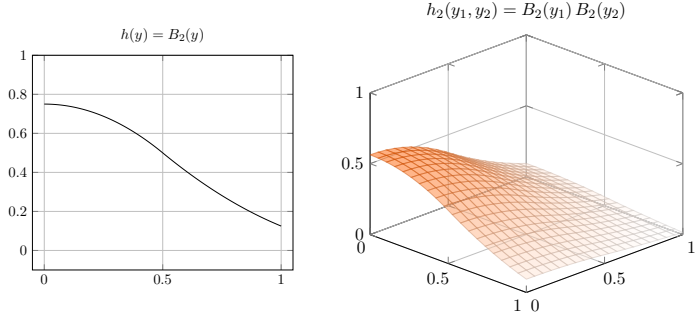

We define a shifted, scaled and dilated B-spline of second order as
\begin{align} \label{def:B2}
	B_2(x) :=
	\begin{cases}
		-x^2 + \frac{3}{4} \quad & \text{ for } 0 \leq x < \frac{1}{2}, \\
		\frac{1}{2}\left( x^2 -3x + \frac{9}{4} \right) \quad & \text{ for } \frac{1}{2} \leq x \leq 1,
	\end{cases}
\end{align}
which we refer to as the \emph{$B_2$-cutoff}, that was also used in \cite{PoVo15,NaPo19}.
It is in $\mathcal{C}^{1}([0,1])$ and depicted in Figure~\ref{fig:B2plot}.
Even though it is only once continuously differentiable, it is also an element in $\mathcal{H}^{\frac{5}{2}-\varepsilon}\left([0,1]\right)$ for any $\varepsilon>0$, which the following arguments show.
It's well-known a second order B-spline is the result of a convolution of three step functions $\chi_{[0,1]}$ (where $\chi$ denotes the indicator function) with themselves, whose respective Fourier coefficients $(\chi_{[0,1]}(\cdot),\mathrm{e}^{2\pi\mathrm{i}k(\cdot)})_{L_2\left([0,1]\right)}$ decay like $|k|^{-1}$ for $k\to\pm\infty$.
Hence, the Fourier coefficients $\hat{h}_k = (B_2,\mathrm{e}^{2\pi\mathrm{i}k(\cdot)})_{L_2\left([0,1]\right)}$ of the $B_2$-cutoff \eqref{def:B2} decay like $|k|^{-3}$ for $k\to\pm\infty$.
Considering a constant weight function $\omega\equiv 1$, the $\|\cdot\|_{\mathcal{H}^{\beta}\left([0,1]\right)}$-norm given in \eqref{def:Sobolev_space_cube} of $B_2$ is finite if
\begin{align*}
	\|B_2\|_{\mathcal{H}^{\beta}\left([0,1]\right)}^2
		= \sum_{k\in\mathbb{Z}} \omega_{\mathrm{hc}}(k)^{2\beta} |\hat h_{k}|^2
		\lesssim \sum_{k\in\mathbb{Z}} \max\{1,|k|\}^{2\beta} \frac{1}{|k|^6} < \infty,
\end{align*}
which is the case for
\begin{align*}
	 |k|^{2\beta-6} \leq k^{-(1+\varepsilon)}
	 \Leftrightarrow \beta \leq \frac{5}{2}-\varepsilon, \quad \varepsilon > 0.
\end{align*}

Next, we approximate the tensored $B_2$-cutoff
\begin{align} \label{def:test_functions}
	h(\mathbf x) = \prod_{\ell=1}^{d} B_2(x_\ell) \in \mathcal{H}^{\frac{5}{2}-\varepsilon}([0,1]^d), \varepsilon>0,
\end{align}
by the approximated Chebyshev, cosine or transformed Fourier partial sums $S_{I}^{\Lambda}h$ given in \eqref{def:approx_cosine_Partsum}, \eqref{def:approx_Cheb_Partsum} and \eqref{def:approx_trafo_Fou_Partsum}.
We study the resulting relative $\ell_2$-and $\ell_{\infty}$-approximation errors
\begin{align} \label{def:ell_p_approxError}
	\varepsilon_{p}^{R}(h)
	:= \frac{\left\| \left( h(\mathbf x_j) - S_{I}^{\Lambda}h(\mathbf x_j) \right)_{j=1}^{R} \right\|_{\ell_{p}}}{\left\| \left( h(\mathbf x_j) \right)_{j=1}^{R}  \right\|_{\ell_{p}}},
	\quad p\in\{2,\infty\},
\end{align}
that are evaluated at $R\in\mathbb{N}$ uniformly distributed points $\mathbf x_j \sim \mathcal{U}([0,1]^d)$.
The approximated coefficients appearing in the approximated partial sums \eqref{def:approx_Cheb_Partsum},\eqref{def:approx_cosine_Partsum} and \eqref{def:approx_trafo_Fou_Partsum} are calculated by solving the corresponding systems \eqref{eq:Chebyshev_system_MatVec}, \eqref{eq:Tent_system_MatVec} or \eqref{eq:Trafo_Exp_system_MatVec}.

\subsection{The numerical results of $\ell_2$-approximation}
Throughout this section we repeatedly use the bold number notation $\mathbf 1 = (1,\ldots,1)^{\top}$ that we already used in the definition of rank-$1$ lattices \eqref{def:R1L} and expressions like $\bm\eta = \mathbf 2$ mean that $\eta_\ell = 2$ for all $\ell\in\{1,\ldots,d\}$.

In \cite{Tem86,ByKaUlVo16,volkmerdiss} we find a broad discussion on the approximation error decay of function in the Sobolev space $\mathcal{H}^{\beta}(\mathbb{T}^d), m\in\mathbb{N}_{0}$.
It was proven that there is a worst case upper error bound of the form 
\begin{align}\label{eq:worstcaseUpperBound}
	\varepsilon_2 = \varepsilon_{2}(N,d)\approx \left\| h - S_{I_{N}^{d}}^{\Lambda}h \right\|_{L_{2}\left(\left[0,1\right]^d\right)} \lesssim N^{-m} (\log N)^{(d-1)/2}.
\end{align}
In \cite{NaPo19} we find conditions on the logarithmic and the error function transformation $\psi(\cdot,\bm\eta)$, given in \eqref{eq:logarithmic_trafo} and \eqref{eq:erf_trafo}, such that a certain degree of smoothness of the given $\mathcal{C}^{m}(\mathbb{T}^d)$-function is preserved under composition with $\psi(\cdot,\bm\eta)$ and the resulting periodized function is at least in $\mathcal H^{\ell}(\mathbb{T}^d), \ell \leq m$ and for each $\ell$ it was calculated how large the parameter $\eta$ has to be chosen.
According to the conditions in \cite[Theorem~4]{NaPo19}, the tensored $B_2$-cutoff in \eqref{def:test_functions} is transformed into a function $f\in\mathcal{H}^{0}(\mathbb{T}^d)$ of the form~\eqref{eq:periodization_idea} for all considered torus-to-cube transformations~$\psi(\cdot,\bm\eta)$ with parameters $1 < \eta_\ell \leq 3$,
and into a function $f\in\mathcal{H}^{1}(\mathbb{T}^d)$ for parameters $\eta_\ell > 3$, $\ell\in\{1,\ldots,d\}$.
While these conditions are independent of the particular considered function $h\in\mathcal{C}^{m}(\mathbb{T}^d)$, they are pretty coarse in the sense of not catching the additional smoothness of functions like the $B_2$-cutoff given in $\eqref{def:test_functions}$ which is an almost $\mathcal{H}^{\frac{5}{2}}([0,1]^d)$-function as we showed earlier.
In numerical tests we showcase that in certain setups the Chebyshev coefficients and the transformed Fourier coefficients will indeed decay faster than the worst case upper bound \eqref{eq:worstcaseUpperBound}.

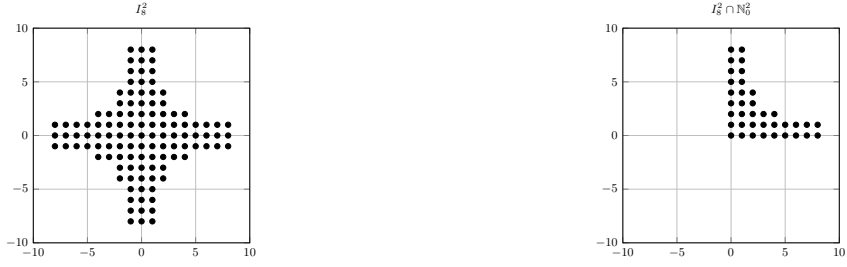
\begin{figure}[tb]
	\begin{minipage}[b]{0.49\linewidth}
		\centering
		\begin{tikzpicture}[scale=0.5]
			\begin{axis}[scatter/classes = { a = {mark=o, draw=black} },
				font=\footnotesize,
				grid = both,
				xmax = 10, xmin = -10, ymax = 10, ymin = -10,
				title = {$I_{8}^{2}$},
unit vector ratio*=1 1 1
				]
				\addplot[scatter ,only marks, mark size=2, scatter src = explicit symbolic]coordinates{
					(-8, -1) (-8, 0) (-8, 1) (-7, -1) (-7, 0) (-7, 1) (-6, -1) (-6, 0) (-6, 1) (-5, -1) (-5, 0) (-5, 1) (-4, -2) (-4, -1) (-4, 0) (-4, 1) (-4, 2) (-3, -2) (-3, -1) (-3, 0) (-3, 1) (-3, 2) (-2, -4) (-2, -3) (-2, -2) (-2, -1) (-2, 0) (-2, 1) (-2, 2) (-2, 3) (-2, 4) (-1, -8) (-1, -7) (-1, -6) (-1, -5) (-1, -4) (-1, -3) (-1, -2) (-1, -1) (-1, 0) (-1, 1) (-1, 2) (-1, 3) (-1, 4) (-1, 5) (-1, 6) (-1, 7) (-1, 8) (0, -8) (0, -7) (0, -6) (0, -5) (0, -4) (0, -3) (0, -2) (0, -1) (0, 0) (0, 1) (0, 2) (0, 3) (0, 4) (0, 5) (0, 6) (0, 7) (0, 8) (1, -8) (1, -7) (1, -6) (1, -5) (1, -4) (1, -3) (1, -2) (1, -1) (1, 0) (1, 1) (1, 2) (1, 3) (1, 4) (1, 5) (1, 6) (1, 7) (1, 8) (2, -4) (2, -3) (2, -2) (2, -1) (2, 0) (2, 1) (2, 2) (2, 3) (2, 4) (3, -2) (3, -1) (3, 0) (3, 1) (3, 2) (4, -2) (4, -1) (4, 0) (4, 1) (4, 2) (5, -1) (5, 0) (5, 1) (6, -1) (6, 0) (6, 1) (7, -1) (7, 0) (7, 1) (8, -1) (8, 0) (8, 1)
				};
			\end{axis}
		\end{tikzpicture}
	\end{minipage}
	\begin{minipage}[b]{0.49\linewidth}
		\centering
		\begin{tikzpicture}[scale=0.5]
			\begin{axis}[scatter/classes = { a = {mark=o, draw=black} },
				font=\footnotesize,
				grid = both,
				xmax = 10, xmin = -10, ymax = 10, ymin = -10,
				title = {$I_{8}^{2}\cap\mathbb{N}_{0}^{2}$},
unit vector ratio*=1 1 1
				]
				\addplot[scatter ,only marks, mark size=2, scatter src = explicit symbolic]coordinates{
					(0, 0) (0, 1) (0, 2) (0, 3) (0, 4) (0, 5) (0, 6) (0, 7) (0, 8) (1, 0) (1, 1) (1, 2) (1, 3) (1, 4) (1, 5) (1, 6) (1, 7) (1, 8) (2, 0) (2, 1) (2, 2) (2, 3) (2, 4) (3, 0) (3, 1) (3, 2) (4, 0) (4, 1) (4, 2)  (5, 0) (5, 1) (6, 0) (6, 1) (7, 0) (7, 1) (8, 0) (8, 1)
				};
			\end{axis}
		\end{tikzpicture}
	\end{minipage}
	\caption{The hyperbolic cross $I_{8}^{2}$ (left) and its first quadrant $I_{8}^{2}\cap\mathbb{N}_{0}^{2}$ (right).
	}
	\label{fig:FreqSets}
\end{figure}

In dimensions $d\in\{1,2,4,7\}$ we compare the discrete $\ell_2$-approximation error $\varepsilon_{2}$, given in \eqref{def:ell_p_approxError}, with $R=1.000.000$ uniformly distributed evaluation points for all of the previously introduced approximation approaches.
We consider frequency sets $I_{N}^{d}$ for all transformed Fourier systems and $I_{N}^{d}\cap\mathbb{N}_{0}^{d}$ for the cosine and Chebyshev systems.
Both frequency sets are illustrated in dimension $d=2$ with $N=8$ in Figure~\ref{fig:FreqSets}.
We use $N\in\{1,\ldots,140\}$ for $d=1$, $N\in\{1,\ldots,80\}$ for $d=2$, $N\in\{1,\ldots,50\}$ for $d=4$ and $N\in\{1,\ldots,30\}$ for $d=7$.

In dimensions $d=1$ and $d=2$ we observe that the approximation errors are significantly better for $\bm\eta=\mathbf 4$ than for $\bm\eta=\mathbf 2$, indicating the increased smoothening effect of both the logarithmic and the error function transformation.
In dimensions $d\in\{4,7\}$, the errors for $\bm\eta=\mathbf 4$ turn out to be worse than for $\bm\eta=\mathbf 2$, which we suspect might be due to the increase of certain constants depending on $\bm\eta$ in the error estimate \eqref{eq:worstcaseUpperBound}.
The Chebyshev approximation turns out to be a solid candidate to approximate the B-spline given in \eqref{def:test_functions}.
In this specific setup, we also checked the error behavior for other parameters $\bm\eta\in \{\mathbf{2.1},\mathbf{2.2},\ldots,\mathbf{3.8},\mathbf{3.9}, \mathbf{4.1},\mathbf{4.2},\ldots\}$.
As it turns out, $\bm\eta = \mathbf{4}$ is the best choice for the logarithmic transformation and for the error function transformation the best choice is $\bm\eta = \mathbf{2.5}$.

However, only the error function transformation is able to match the approximation error of the Chebyshev approximation, which also shows when we investigate and compare the error decay rates of $\varepsilon_2^R(h)$ that were numerically observed for the univariate case $d=1$.
In this specific setup, $h$ is still the continuous second-order B-spline given in \eqref{def:test_functions} that is an element of 
$\mathcal{H}^{\frac{5}{2}-\varepsilon}\left([0,1]^d\right)$.
Hence, we expect to obtain an error decay at most $\varepsilon_2^R(h)\lesssim N^{-\frac{5}{2}+\varepsilon}$ for any $\varepsilon > 0$ and increasing values of $N$ when approximating $h$ with respect to any transformed Fourier system.
We achieve these decay rates numerically with the Chebyshev system and with the transformed Fourier system when considering the logarithmic transformation with $\eta \in\{2.5,4\}$.
Interestingly, the decay rates of the cosine system remain at $N^{-1.5}$.
In comparison, the logarithmically transformed Fourier system with $\eta = 2$ loses half an order, which is slightly improved for $\eta = 4$.
In total we observe that some transformed Fourier systems are able to achieve the same decay rates as the Chebyshev system, when we use parameterized torus-to-cube transformations $\psi(\cdot,\eta)$ and pick an appropriate parameter $\eta\in\mathbb{R}_+$.
The results are summarized in Table~\ref{table:CompErrorsB2}.
\begin{table}
\centering
\begin{tabular}{|l|p{1.5cm}|}
	\hline
transformation & $\varepsilon_2^R(h)$  \\
	\hline
	\eqref{def:cosine_system} cosine system 										& $N^{-1.5}$ 	\\
	\eqref{def:T_k} Chebyshev system 												& $N^{-2.45}$ 	 \\
	\eqref{eq:logarithmic_trafo} log transf. Fourier, $\eta = 2$				& $N^{-1}$ 		 \\
	\eqref{eq:logarithmic_trafo} log transf. Fourier, $\eta = 4$				& $N^{-2.25}$ 	 \\
	\eqref{eq:erf_trafo} error fct. transf. Fourier, $\eta = 2$ 	& $N^{-1.9}$ 	 \\
	\eqref{eq:erf_trafo} error fct. transf. Fourier, $\eta = 2.5$ 	& $N^{-2.5}$ 	 \\
	\eqref{eq:erf_trafo} error fct. transf. Fourier, $\eta = 4$ 	& $N^{-2.5}$ 	 \\
	\hline
\end{tabular}
\caption{
	The observed decay rates of the discrete approximation error $\varepsilon_{2}^{R}(h)$ as given in \eqref{def:ell_p_approxError} 
	when $h$ is the univariate $B_2$-cutoff as defined in \eqref{def:B2}.
}
\label{table:CompErrorsB2}
\end{table}

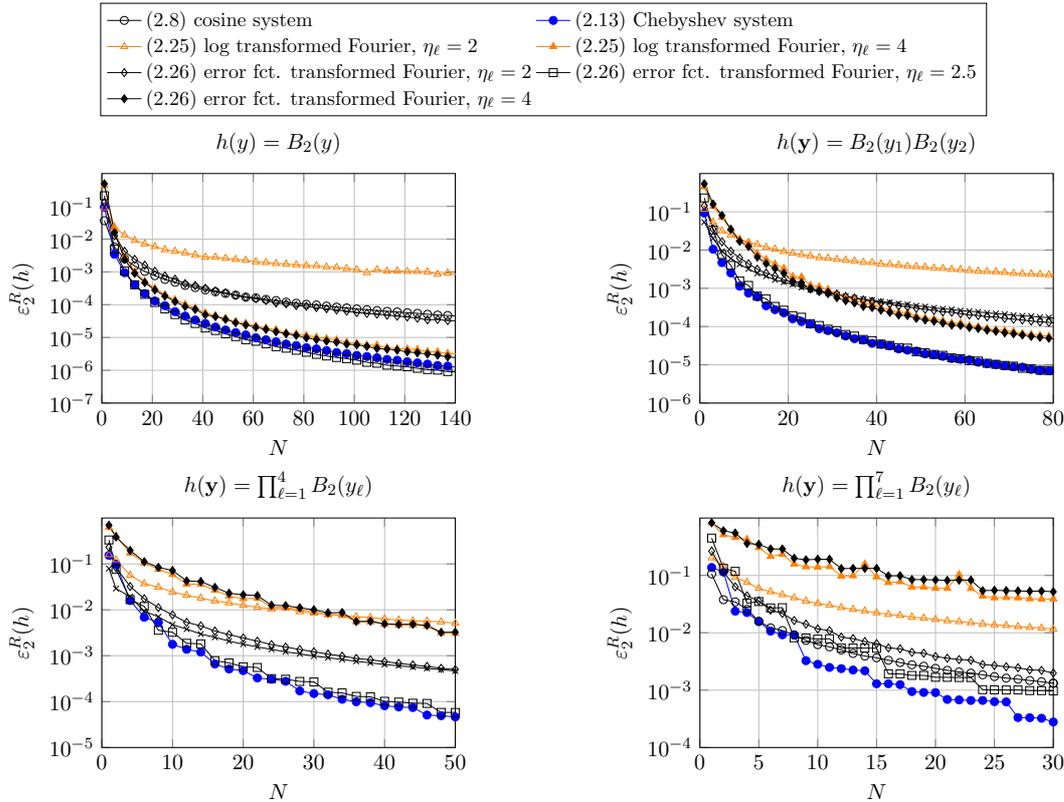
\begin{figure}[tb]
	\begin{minipage}[b]{0.5\linewidth}
		\begin{tikzpicture}[scale=0.75]
		\begin{axis}[
			width=1.0\textwidth, height=0.725\textwidth,
			ymode = log,
			enlargelimits=false,
			xmin=0, xmax=140, ymin=1e-7, ymax=1e-0,
			ytick={1e-1,1e-2,1e-3,1e-4,1e-5,1e-6,1e-7,1e-8,1e-9,1e-10,1e-11},
			grid=major, 
			title ={$h(y) = B_2(y)$},
			xlabel={$N$}, 
			ylabel={$\varepsilon_2^R(h)$},
			legend style={at={(1.25,1.25)}, anchor=south,legend columns=2,legend cell align=left, font=\small},
			xminorticks=false,
			yminorticks=false
			]
			\addplot[black, line width = 0.5pt, mark=o] coordinates {
				(1, 3.6517e-02) (5, 5.7119e-03) (9, 2.5086e-03) (13, 1.5151e-03) (17, 1.0392e-03) (21, 7.6610e-04) (25, 5.8732e-04) (29, 4.8287e-04) (33, 3.8855e-04) (37, 3.3260e-04) (41, 2.8246e-04) (45, 2.4677e-04) (49, 2.2047e-04) (53, 1.9629e-04) (57, 1.7049e-04) (61, 1.6188e-04) (65, 1.4509e-04) (69, 1.3334e-04) (73, 1.1978e-04) (77, 1.1074e-04) (81, 1.0283e-04) (85, 9.8015e-05) (89, 9.0795e-05) (93, 8.4340e-05) (97, 7.9616e-05) (101, 7.3969e-05) (105, 7.0398e-05) (109, 6.6159e-05) (113, 6.2830e-05) (117, 6.0147e-05) (121, 5.6819e-05) (125, 5.4457e-05) (129, 5.1108e-05) (133, 4.8313e-05) (137, 4.6881e-05) (141, 4.5286e-05) (145, 4.2355e-05) (149, 4.1271e-05) (153, 4.0276e-05) (157, 3.7762e-05) (161, 3.7606e-05) (165, 3.5756e-05) (169, 3.1020e-05) (173, 3.3382e-05) (177, 3.2476e-05) (181, 3.0821e-05) (185, 3.0547e-05) (189, 2.8017e-05) (193, 2.8871e-05) (197, 2.6805e-05) (201, 2.6437e-05)
			};
			\addlegendentry{\eqref{def:cosine_system} cosine system};
			\addplot[blue, line width = 0.5pt, mark=*] coordinates {
				(1, 9.4790e-02) (5, 3.4471e-03) (9, 9.4422e-04) (13, 4.0451e-04) (17, 2.1851e-04) (21, 1.2974e-04) (25, 9.3300e-05) (29, 6.1301e-05) (33, 4.4891e-05) (37, 3.4109e-05) (41, 2.6610e-05) (45, 2.0807e-05) (49, 1.6940e-05) (53, 1.4090e-05) (57, 1.1962e-05) (61, 9.7963e-06) (65, 8.5162e-06) (69, 7.3513e-06) (73, 6.3251e-06) (77, 5.6904e-06) (81, 4.8787e-06) (85, 4.4831e-06) (89, 3.9654e-06) (93, 3.4756e-06) (97, 3.1946e-06) (101, 2.8388e-06) (105, 2.6279e-06) (109, 2.4078e-06) (113, 2.1460e-06) (117, 2.0150e-06) (121, 1.8179e-06) (125, 1.6796e-06) (129, 1.5435e-06) (133, 1.3987e-06) (137, 1.3397e-06) (141, 1.2291e-06) (145, 1.1459e-06) (149, 1.1042e-06) (153, 1.0553e-06) (157, 9.8583e-07) (161, 9.2242e-07) (165, 8.4641e-07) (169, 7.8489e-07) (173, 7.6726e-07) (177, 7.0727e-07) (181, 6.6456e-07) (185, 6.5669e-07) (189, 5.9454e-07) (193, 5.7439e-07) (197, 5.5338e-07) (201, 5.2465e-07)
			};
			\addlegendentry{\eqref{def:T_k} Chebyshev system};
			\addplot[mark=triangle, orange, line width = 0.5pt] coordinates {
				(1, 9.0042e-02)  (5, 2.3151e-02)  (9, 1.3263e-02)  (13, 9.1829e-03)  (17, 7.1609e-03)  (21, 5.9348e-03)  (25, 4.7652e-03)  (29, 4.2930e-03)  (33, 3.8177e-03)  (37, 3.3127e-03)  (41, 2.8064e-03)  (45, 2.8147e-03)  (49, 2.5878e-03)  (53, 2.4106e-03)  (57, 2.2632e-03)  (61, 2.0230e-03)  (65, 1.9386e-03)  (69, 1.8039e-03)  (73, 1.7390e-03)  (77, 1.6439e-03)  (81, 1.5588e-03)  (85, 1.5390e-03)  (89, 1.4208e-03)  (93, 1.3668e-03)  (97, 1.2655e-03)  (101, 1.1694e-03)  (105, 9.4735e-04)  (109, 1.1725e-03)  (113, 1.0527e-03)  (117, 1.0679e-03)  (121, 1.0366e-03)  (125, 1.0368e-03)  (129, 1.0000e-03)  (133, 8.4715e-04)  (137, 9.2816e-04)  (141, 8.2719e-04)  (145, 8.4158e-04)  (149, 7.3399e-04)  (153, 8.3494e-04)  (157, 7.5223e-04)  (161, 7.3858e-04)  (165, 7.6447e-04)  (169, 7.4776e-04)  (173, 6.8208e-04)  (177, 6.7957e-04)  (181, 6.2375e-04)  (185, 6.8833e-04)  (189, 5.8919e-04)  (193, 5.7743e-04)  (197, 5.8959e-04)  (201, 6.1573e-04)
			};
			\addlegendentry{\eqref{eq:logarithmic_trafo} log transformed Fourier, $\eta_\ell = 2$};
			\addplot[mark=triangle*, orange, line width = 0.5pt] coordinates {
				(1, 4.2738e-01)  (5, 2.0151e-02)  (9, 2.8848e-03)  (13, 1.0199e-03)  (17, 5.1724e-04)  (21, 3.1734e-04)  (25, 1.9958e-04)  (29, 1.4062e-04)  (33, 1.0269e-04)  (37, 7.8701e-05)  (41, 5.6510e-05)  (45, 4.6411e-05)  (49, 3.8228e-05)  (53, 3.0235e-05)  (57, 2.7106e-05)  (61, 2.2224e-05)  (65, 1.9415e-05)  (69, 1.7284e-05)  (73, 1.3924e-05)  (77, 1.2838e-05)  (81, 1.1749e-05)  (85, 1.0267e-05)  (89, 9.3010e-06)  (93, 8.3665e-06)  (97, 7.1821e-06)  (101, 6.6960e-06)  (105, 5.8256e-06)  (109, 5.8123e-06)  (113, 4.8684e-06)  (117, 4.5929e-06)  (121, 4.4565e-06)  (125, 4.2560e-06)  (129, 3.7125e-06)  (133, 3.6551e-06)  (137, 3.2424e-06)  (141, 3.1918e-06)  (145, 2.9328e-06)  (149, 2.7470e-06)  (153, 2.5952e-06)  (157, 2.1667e-06)  (161, 2.3560e-06)  (165, 2.2082e-06)  (169, 2.0119e-06)  (173, 1.9466e-06)  (177, 1.9265e-06)  (181, 1.7641e-06)  (185, 1.6705e-06)  (189, 1.6100e-06)  (193, 1.4904e-06)  (197, 1.4384e-06)  (201, 1.4222e-06)  
			};
			\addlegendentry{\eqref{eq:logarithmic_trafo} log transformed Fourier, $\eta_\ell = 4$}; 
			\addplot[mark=diamond, line width = 0.5pt] coordinates {
				(1, 1.3004e-01)  (5, 1.2821e-02)  (9, 4.2360e-03)  (13, 2.4107e-03)  (17, 1.6028e-03)  (21, 1.0462e-03)  (25, 7.5270e-04)  (29, 5.9708e-04)  (33, 4.3850e-04)  (37, 3.8394e-04)  (41, 3.1757e-04)  (45, 2.6240e-04)  (49, 2.3013e-04)  (53, 1.9200e-04)  (57, 1.6854e-04)  (61, 1.5190e-04)  (65, 1.3460e-04)  (69, 1.1305e-04)  (73, 1.0690e-04)  (77, 9.4208e-05)  (81, 8.9317e-05)  (85, 8.1696e-05)  (89, 7.5550e-05)  (93, 6.9102e-05)  (97, 6.3473e-05)  (101, 5.7696e-05)  (105, 5.4191e-05)  (109, 5.1704e-05)  (113, 4.8144e-05)  (117, 4.4773e-05)  (121, 4.2276e-05)  (125, 3.8786e-05)  (129, 3.4619e-05)  (133, 3.5420e-05)  (137, 3.2618e-05)  (141, 3.2399e-05)  (145, 3.0972e-05)  (149, 2.7569e-05)  (153, 2.6688e-05)  (157, 2.5727e-05)  (161, 2.3757e-05)  (165, 2.4147e-05)  (169, 2.3202e-05)  (173, 2.2193e-05)  (177, 1.9573e-05)  (181, 1.9234e-05)  (185, 1.9109e-05)  (189, 1.7976e-05)  (193, 1.6955e-05)  (197, 1.4908e-05)  (201, 1.6738e-05) 
			};
			\addlegendentry{\eqref{eq:erf_trafo} error fct. transformed Fourier, $\eta_\ell = 2$}; 
			\addplot[mark=square, line width = 0.5pt] coordinates {
				(1, 2.1051e-01)  (5, 4.8835e-03)  (9, 1.0647e-03)  (13, 4.0451e-04)  (17, 2.0643e-04)  (21, 1.1650e-04)  (25, 7.4148e-05)  (29, 4.9968e-05)  (33, 3.4270e-05)  (37, 2.6695e-05)  (41, 1.9397e-05)  (45, 1.6284e-05)  (49, 1.2908e-05)  (53, 1.0653e-05)  (57, 8.5309e-06)  (61, 7.6579e-06)  (65, 6.4519e-06)  (69, 5.2414e-06)  (73, 4.7115e-06)  (77, 4.0764e-06)  (81, 3.5418e-06)  (85, 3.0547e-06)  (89, 2.7272e-06)  (93, 2.4505e-06)  (97, 2.2608e-06)  (101, 2.0320e-06)  (105, 1.7501e-06)  (109, 1.6396e-06)  (113, 1.4293e-06)  (117, 1.3464e-06)  (121, 1.2749e-06)  (125, 1.1502e-06)  (129, 1.0518e-06)  (133, 1.0115e-06)  (137, 9.1424e-07)  (141, 8.8666e-07)  (145, 7.7990e-07)  (149, 7.2082e-07)  (153, 6.7309e-07)  (157, 6.5369e-07)  (161, 6.0420e-07)  (165, 5.9003e-07)  (169, 5.2693e-07)  (173, 4.9742e-07)  (177, 4.8428e-07)  (181, 4.5394e-07)  (185, 4.3191e-07)  (189, 4.1982e-07)  (193, 3.6179e-07)  (197, 3.7131e-07)  (201, 3.4487e-07)
			};
			\addlegendentry{\eqref{eq:erf_trafo} error fct. transformed Fourier, $\eta_\ell = 2.5$}; 
			\addplot[mark=diamond*, line width = 0.5pt] coordinates {
				(1, 4.8733e-01)  (5, 1.5822e-02)  (9, 2.3978e-03)  (13, 9.0415e-04)  (17, 4.7397e-04)  (21, 2.8068e-04)  (25, 1.8365e-04)  (29, 1.2449e-04)  (33, 8.7525e-05)  (37, 7.1781e-05)  (41, 4.9091e-05)  (45, 4.3412e-05)  (49, 3.5641e-05)  (53, 2.8230e-05)  (57, 2.4103e-05)  (61, 2.0443e-05)  (65, 1.7643e-05)  (69, 1.4953e-05)  (73, 1.3073e-05)  (77, 1.1633e-05)  (81, 9.7939e-06)  (85, 8.6098e-06)  (89, 7.7199e-06)  (93, 7.1941e-06)  (97, 6.2519e-06)  (101, 5.9670e-06)  (105, 5.3453e-06)  (109, 4.6117e-06)  (113, 4.3860e-06)  (117, 4.1110e-06)  (121, 3.6502e-06)  (125, 3.3933e-06)  (129, 3.1271e-06)  (133, 2.8851e-06)  (137, 2.5508e-06)  (141, 2.4139e-06)  (145, 2.4338e-06)  (149, 2.0891e-06)  (153, 2.0775e-06)  (157, 1.9916e-06)  (161, 1.8550e-06)  (165, 1.7004e-06)  (169, 1.5054e-06)  (173, 1.5686e-06)  (177, 1.3050e-06)  (181, 1.3334e-06)  (185, 1.3229e-06)  (189, 1.1631e-06)  (193, 1.1732e-06)  (197, 1.1183e-06)  (201, 1.0066e-06)  
			};
			\addlegendentry{\eqref{eq:erf_trafo} error fct. transformed Fourier, $\eta_\ell = 4$};
			\end{axis}
		\end{tikzpicture}
	\end{minipage}
	\begin{minipage}[b]{0.5\linewidth}
		\begin{tikzpicture}[scale=0.75]
		\begin{axis}[
			width=1.0\textwidth, height=0.725\textwidth,
			ymode = log,
			enlargelimits=false,
			xmin=0, xmax=80, ymin=1e-6, ymax=1e-0,
			ytick={1e-1,1e-2,1e-3,1e-4,1e-5,1e-6,1e-7,1e-8,1e-9,1e-10,1e-11},
			grid=major, 
			title ={$h(\mathbf y) = B_2(y_1)B_2(y_2)$},
			xlabel={$N$}, 
			ylabel={$\varepsilon_2^R(h)$},
xminorticks=false,
			yminorticks=false
			]
			\addplot[blue, line width = 0.5pt, mark=*] coordinates {
				(1, 9.2684e-02)  (3, 1.0512e-02)  (5, 4.6930e-03)  (7, 2.5294e-03)  (9, 1.1349e-03)  (11, 7.4937e-04)  (13, 5.9793e-04)  (15, 3.4845e-04)  (17, 2.7607e-04)  (19, 2.2828e-04)  (21, 1.5918e-04)  (23, 1.3496e-04)  (25, 1.1558e-04)  (27, 8.7555e-05)  (29, 7.8352e-05)  (31, 6.8274e-05)  (33, 5.6692e-05)  (35, 4.8118e-05)  (37, 4.4481e-05)  (39, 3.6319e-05)  (41, 3.4342e-05)  (43, 3.1444e-05)  (45, 2.6016e-05)  (47, 2.4126e-05)  (49, 2.2737e-05)  (51, 1.9343e-05)  (53, 1.8273e-05)  (55, 1.6649e-05)  (57, 1.4683e-05)  (59, 1.3954e-05)  (61, 1.3427e-05)  (63, 1.1685e-05)  (65, 1.0755e-05)  (67, 1.0230e-05)  (69, 9.4127e-06)  (71, 8.8554e-06)  (73, 8.6633e-06)  (75, 7.6137e-06)  (77, 7.1664e-06)  (79, 6.9390e-06)  (81, 6.3346e-06)  
			};
\addplot[black, line width = 0.5pt, mark=x] coordinates {
				(1, 5.3674e-02)  (3, 2.0043e-02)  (5, 9.2807e-03)  (7, 6.0622e-03)  (9, 4.2340e-03)  (11, 3.2133e-03)  (13, 2.4448e-03)  (15, 2.0737e-03)  (17, 1.6559e-03)  (19, 1.4477e-03)  (21, 1.2639e-03)  (23, 1.0565e-03)  (25, 9.1370e-04)  (27, 8.2773e-04)  (29, 7.5289e-04)  (31, 6.7976e-04)  (33, 6.2343e-04)  (35, 5.6445e-04)  (37, 5.1874e-04)  (39, 4.9388e-04)  (41, 4.4740e-04)  (43, 4.3059e-04)  (45, 3.9622e-04)  (47, 3.7488e-04)  (49, 3.4751e-04)  (51, 3.3153e-04)  (53, 3.0728e-04)  (55, 2.9022e-04)  (57, 2.7954e-04)  (59, 2.6088e-04)  (61, 2.5139e-04)  (63, 2.4196e-04)  (65, 2.2691e-04)  (67, 2.1916e-04)  (69, 2.1041e-04)  (71, 1.9815e-04)  (73, 1.9002e-04)  (75, 1.8641e-04)  (77, 1.7561e-04)  (79, 1.7168e-04)  (81, 1.6236e-04) 
			};
\addplot[mark=triangle, orange, line width = 0.5pt] coordinates {
				(1, 1.1619e-01)  (3, 5.3715e-02)  (5, 3.1781e-02)  (7, 2.4788e-02)  (9, 1.9237e-02)  (11, 1.5887e-02)  (13, 1.4056e-02)  (15, 1.1789e-02)  (17, 1.0576e-02)  (19, 9.1362e-03)  (21, 8.2696e-03)  (23, 7.6162e-03)  (25, 7.0426e-03)  (27, 6.3469e-03)  (29, 6.1456e-03)  (31, 5.7418e-03)  (33, 5.4687e-03)  (35, 5.1750e-03)  (37, 4.7908e-03)  (39, 4.7164e-03)  (41, 4.2964e-03)  (43, 4.1995e-03)  (45, 4.0124e-03)  (47, 3.7781e-03)  (49, 3.6953e-03)  (51, 3.5747e-03)  (53, 3.3688e-03)  (55, 3.3186e-03)  (57, 3.1281e-03)  (59, 3.0678e-03)  (61, 2.9240e-03)  (63, 2.8398e-03)  (65, 2.7516e-03)  (67, 2.7143e-03)  (69, 2.6045e-03)  (71, 2.4771e-03)  (73, 2.4129e-03)  (75, 2.3855e-03)  (77, 2.3163e-03)  (79, 2.2580e-03)  (81, 2.1794e-03)  
			};
\addplot[mark=triangle*, orange, line width = 0.5pt] coordinates {
				(1, 4.4948e-01)  (3, 1.3840e-01)  (5, 7.1451e-02)  (7, 3.5606e-02)  (9, 1.9152e-02)  (11, 1.4700e-02)  (13, 8.1160e-03)  (15, 5.6728e-03)  (17, 4.2955e-03)  (19, 3.3339e-03)  (21, 2.0280e-03)  (23, 1.9740e-03)  (25, 1.2004e-03)  (27, 1.1022e-03)  (29, 9.3131e-04)  (31, 7.1941e-04)  (33, 5.7449e-04)  (35, 5.0929e-04)  (37, 4.0786e-04)  (39, 3.8065e-04)  (41, 3.1198e-04)  (43, 2.7304e-04)  (45, 2.3557e-04)  (47, 2.3413e-04)  (49, 1.8546e-04)  (51, 1.6669e-04)  (53, 1.5888e-04)  (55, 1.4059e-04)  (57, 1.2392e-04)  (59, 1.2303e-04)  (61, 1.0491e-04)  (63, 9.9283e-05)  (65, 8.8650e-05)  (67, 8.2564e-05)  (69, 7.9110e-05)  (71, 7.4004e-05)  (73, 6.4802e-05)  (75, 6.2299e-05)  (77, 5.8732e-05)  (79, 5.6039e-05)  (81, 4.9739e-05)  
			};
\addplot[mark=diamond, line width = 0.5pt] coordinates {
				(1, 1.4338e-01)  (3, 3.5207e-02)  (5, 1.6422e-02)  (7, 9.3125e-03)  (9, 6.1806e-03)  (11, 4.3697e-03)  (13, 3.1550e-03)  (15, 2.4943e-03)  (17, 2.1029e-03)  (19, 1.6609e-03)  (21, 1.4395e-03)  (23, 1.1936e-03)  (25, 9.8672e-04)  (27, 8.7066e-04)  (29, 7.8550e-04)  (31, 7.0147e-04)  (33, 6.1285e-04)  (35, 5.5639e-04)  (37, 4.9589e-04)  (39, 4.5611e-04)  (41, 4.1875e-04)  (43, 3.8770e-04)  (45, 3.4849e-04)  (47, 3.2805e-04)  (49, 3.0070e-04)  (51, 2.8274e-04)  (53, 2.6480e-04)  (55, 2.4267e-04)  (57, 2.3707e-04)  (59, 2.1697e-04)  (61, 2.0134e-04)  (63, 1.9456e-04)  (65, 1.7934e-04)  (67, 1.7779e-04)  (69, 1.6315e-04)  (71, 1.5547e-04)  (73, 1.4866e-04)  (75, 1.3939e-04)  (77, 1.3241e-04)  (79, 1.2955e-04)  (81, 1.2022e-04)
			};
\addplot[mark=square, line width = 0.5pt] coordinates {
				(1, 2.3034e-01)  (3, 3.3322e-02)  (5, 7.5782e-03)  (7, 3.9727e-03)  (9, 1.5483e-03)  (11, 1.1888e-03)  (13, 6.2598e-04)  (15, 5.3032e-04)  (17, 3.3602e-04)  (19, 2.6782e-04)  (21, 1.9230e-04)  (23, 1.7508e-04)  (25, 1.1893e-04)  (27, 1.0653e-04)  (29, 8.5259e-05)  (31, 7.8052e-05)  (33, 6.2030e-05)  (35, 5.6354e-05)  (37, 4.5913e-05)  (39, 4.2792e-05)  (41, 3.4833e-05)  (43, 3.2164e-05)  (45, 2.8169e-05)  (47, 2.7288e-05)  (49, 2.1996e-05)  (51, 2.0662e-05)  (53, 1.8682e-05)  (55, 1.7665e-05)  (57, 1.5165e-05)  (59, 1.4477e-05)  (61, 1.2857e-05)  (63, 1.2513e-05)  (65, 1.0958e-05)  (67, 1.0280e-05)  (69, 9.6104e-06)  (71, 9.2249e-06)  (73, 8.1174e-06)  (75, 7.8499e-06)  (77, 7.1748e-06)  (79, 7.1229e-06)  (81, 6.2092e-06) 
			};
\addplot[mark=diamond*,line width = 0.5pt] coordinates {
				(1, 5.4031e-01)  (3, 1.5886e-01)  (5, 8.0031e-02)  (7, 3.4670e-02)  (9, 1.6954e-02)  (11, 1.2383e-02)  (13, 6.5300e-03)  (15, 4.4230e-03)  (17, 3.2928e-03)  (19, 2.4691e-03)  (21, 1.5324e-03)  (23, 1.3996e-03)  (25, 9.5604e-04)  (27, 8.2026e-04)  (29, 7.2677e-04)  (31, 5.6950e-04)  (33, 4.7442e-04)  (35, 4.1408e-04)  (37, 3.3720e-04)  (39, 3.1426e-04)  (41, 2.6933e-04)  (43, 2.3440e-04)  (45, 2.0004e-04)  (47, 1.9394e-04)  (49, 1.6333e-04)  (51, 1.4910e-04)  (53, 1.3791e-04)  (55, 1.2416e-04)  (57, 1.0935e-04)  (59, 1.0702e-04)  (61, 9.4351e-05)  (63, 8.8958e-05)  (65, 8.0525e-05)  (67, 7.5246e-05)  (69, 6.8517e-05)  (71, 6.4691e-05)  (73, 6.0032e-05)  (75, 5.5203e-05)  (77, 5.2177e-05)  (79, 4.9477e-05)  (81, 4.4935e-05)      
			};
\end{axis}
		\end{tikzpicture}
	\end{minipage}

	\begin{minipage}[b]{0.5\linewidth}
		\begin{tikzpicture}[scale=0.75]
		\begin{axis}[
			width=1.0\textwidth, height=0.725\textwidth,
			ymode = log,
			enlargelimits=false,
			xmin=0, xmax=50, ymin=1e-5, ymax=1e-0,
			ytick={1e-1,1e-2,1e-3,1e-4,1e-5,1e-6,1e-7,1e-8,1e-9,1e-10,1e-11},
			grid=major, 
			title ={$h(\mathbf y) = \prod_{\ell=1}^{4} B_2(y_\ell)$},
			xlabel={$N$}, 
			ylabel={$\varepsilon_2^R(h)$},
xminorticks=false,
			yminorticks=false
			]
			\addplot[blue, line width = 0.5pt, mark=*] coordinates {
				(1, 1.5336e-01) (2, 9.8351e-02) (4, 1.5902e-02) (6, 6.8732e-03)(8, 5.3262e-03)  (10, 1.7716e-03) (12, 1.3783e-03)(14, 1.2004e-03)(16, 6.5982e-04)  (18, 5.1376e-04)(20, 4.7307e-04)  (22, 3.2802e-04)  (24, 3.0705e-04) (26, 2.7695e-04) (28, 1.7079e-04) (30, 1.4981e-04)  (32, 1.4070e-04) (34, 1.1242e-04)(36, 9.9294e-05) (38, 9.4916e-05)  (40, 8.1103e-05) (42, 7.6393e-05) (44, 7.4186e-05)  (46, 5.0976e-05)(48, 4.9233e-05)  (50, 4.6394e-05)
			};
\addplot[black, line width = 0.5pt, mark=x] coordinates {
				(1, 7.9412e-02) (2, 2.8975e-02) (4, 1.9354e-02) (6, 9.2842e-03) (8, 7.0494e-03) (10, 4.7395e-03) (12, 3.7974e-03) (14, 2.9644e-03)  (16, 2.5138e-03)(18, 2.0782e-03)  (20, 1.7826e-03)  (22, 1.5503e-03)(24, 1.3749e-03)  (26, 1.2130e-03)(28, 1.0994e-03) (30, 9.8050e-04) (32, 9.0444e-04)  (34, 8.2639e-04)  (36, 7.5536e-04) (38, 6.9935e-04) (40, 6.4148e-04) (42, 6.0290e-04) (44, 5.6498e-04)  (46, 5.2656e-04) (48, 4.9646e-04) (50, 4.6502e-04) 
			};
\addplot[mark=triangle, orange, line width = 0.5pt] coordinates {
				(1, 1.5872e-01)  (2, 1.2316e-01)   (4, 5.6765e-02)   (6, 3.7063e-02)   (8, 3.1256e-02)   (10, 2.4349e-02) (12, 2.0771e-02)    (14, 1.8004e-02)   (16, 1.5621e-02)    (18, 1.4168e-02)   (20, 1.2505e-02)   (22, 1.1337e-02)   (24, 1.0438e-02)   (26, 1.0897e-02)  (28, 9.0727e-03)   (30, 8.4499e-03)  (32, 8.0365e-03)  (34, 7.4951e-03)  (36, 7.0985e-03)  (38, 6.7971e-03)   (40, 6.3793e-03)  (42, 5.9824e-03)   (44, 5.7466e-03)    (46, 5.6734e-03)  (48, 5.4605e-03)  (50, 5.1371e-03)  
			};
\addplot[mark=triangle*, orange, line width = 0.5pt] coordinates {
				(1, 6.2891e-01)  (2, 4.0508e-01)   (4, 1.7173e-01)  (6, 1.0535e-01)   (8, 7.6218e-02)  (10, 5.8807e-02)   (12, 3.6171e-02)(14, 3.5804e-02)   (16, 2.6587e-02)  (18, 1.9749e-02)   (20, 1.7578e-02) (22, 1.7730e-02)  (24, 1.1407e-02)  (26, 1.1196e-02)  (28, 1.0296e-02)  (30, 8.9847e-03)     (32, 7.7894e-03)  (34, 7.7688e-03)   (36, 5.4800e-03)  (38, 5.4497e-03)  (40, 4.8209e-03)   (42, 4.6957e-03)    (44, 4.6478e-03)  (46, 4.4236e-03)  (48, 3.1187e-03)   (50, 3.0647e-03)   
			};
\addplot[mark=diamond, line width = 0.5pt] coordinates {
				(1, 2.3046e-01)  (2, 8.3284e-02)  (4, 3.2228e-02)  (6, 1.7394e-02)  (8, 1.1156e-02)   (10, 7.7491e-03)   (12, 5.3322e-03)   (14, 4.4658e-03) (16, 3.4436e-03) (18, 2.8761e-03)  (20, 2.4449e-03)  (22, 2.0734e-03)   (24, 1.7495e-03) (26, 1.5451e-03)  (28, 1.3657e-03)   (30, 1.2168e-03)  (32, 1.0842e-03)  (34, 9.7990e-04)   (36, 8.7205e-04)   (38, 8.1482e-04)  (40, 7.4512e-04)  (42, 6.8148e-04)   (44, 6.2854e-04)  (46, 5.8579e-04)  (48, 5.3170e-04) (50, 4.9987e-04) 
			};
\addplot[mark=square, line width = 0.5pt] coordinates {
				(1, 3.3233e-01)  (2, 7.7052e-02)  (4, 1.6704e-02)   (6, 1.1853e-02)   (8, 3.5784e-03) (10, 3.2412e-03)  (12, 1.8639e-03)  (14, 1.8104e-03)   (16, 7.8057e-04)  (18, 6.9583e-04)   (20, 5.7324e-04)  (22, 5.5360e-04)  (24, 3.1729e-04)  (26, 3.0758e-04)(28, 2.7508e-04)  (30, 2.6846e-04)(32, 1.6035e-04)    (34, 1.5473e-04) (36, 1.2986e-04)  (38, 1.2694e-04)  (40, 1.0159e-04)  (42, 9.7783e-05)   (44, 9.2422e-05)   (46, 9.1211e-05)  (48, 5.8204e-05) (50, 5.8039e-05)
			};
\addplot[mark=diamond*, line width = 0.5pt] coordinates {
				(1, 6.9610e-01)  (2, 3.8593e-01)  (4, 1.9677e-01)  (6, 1.1157e-01)  (8, 8.4091e-02)  (10, 7.1837e-02)  (12, 4.2603e-02)  (14, 4.1267e-02)  (16, 3.0764e-02)  (18, 2.2929e-02)  (20, 2.1200e-02)  (22, 2.0403e-02)  (24, 1.2572e-02)  (26, 1.2312e-02)  (28, 1.1101e-02)  (30, 9.9975e-03)  (32, 8.5922e-03)  (34, 8.7264e-03)  (36, 5.5593e-03)  (38, 5.6492e-03)  (40, 4.9617e-03)  (42, 4.8103e-03)  (44, 4.8655e-03)  (46, 4.5832e-03)  (48, 3.2273e-03)  (50, 3.2621e-03)  
			};
\end{axis}
		\end{tikzpicture}
	\end{minipage}
	\begin{minipage}[b]{0.5\linewidth}
		\begin{tikzpicture}[scale=0.75]
		\begin{axis}[
			width=1.0\textwidth, height=0.725\textwidth,
			ymode = log,
			enlargelimits=false,
			xmin=0, xmax=30, ymin=1e-4, ymax=1e-0,
			ytick={1e-1,1e-2,1e-3,1e-4,1e-5,1e-6,1e-7,1e-8,1e-9,1e-10,1e-11},
			grid=major, 
			title ={$h(\mathbf y) = \prod_{\ell=1}^{7} B_2(y_\ell)$},
			xlabel={$N$}, 
			ylabel={$\varepsilon_2^R(h)$},
xminorticks=false,
			yminorticks=false
			]
			\addplot[blue, line width = 0.5pt, mark=*] coordinates {
				(1, 1.3737e-01) (2, 1.1497e-01) (3, 2.3769e-02) (4, 2.2477e-02) (5, 1.5694e-02) (6, 1.0682e-02) (7, 9.2206e-03) (8, 9.1945e-03) (9, 3.2651e-03) (10, 2.8229e-03) (11, 2.4903e-03) (12, 2.3715e-03) (13, 2.2416e-03) (14, 2.1740e-03) (15, 1.2923e-03) (16, 1.2956e-03) (17, 1.2468e-03) (18, 9.4465e-04) (19, 9.1150e-04) (20, 9.0607e-04) (21, 6.8875e-04) (22, 6.7624e-04) (23, 6.6166e-04) (24, 6.6171e-04) (25, 6.3037e-04) (26, 6.2425e-04) (27, 3.3278e-04) (28, 3.3114e-04) (29, 3.2374e-04) (30, 2.7871e-04) (31, 2.7296e-04) (32, 2.7284e-04) (33, 2.3168e-04) (34, 2.2919e-04) (35, 2.1114e-04) (36, 2.0550e-04) (37, 2.0288e-04) (38, 2.0114e-04) (39, 1.8168e-04) (40, 1.8170e-04)
			};
\addplot[black, line width = 0.5pt, mark=o] coordinates {
				(1, 1.0586e-01) (2, 3.7701e-02) (3, 3.4534e-02) (4, 2.5648e-02) (5, 1.5857e-02) (6, 1.2170e-02) (7, 1.0938e-02) (8, 9.1452e-03) (9, 7.3734e-03) (10, 6.1984e-03) (11, 5.6759e-03) (12, 5.0166e-03) (13, 4.3556e-03) (14, 3.9290e-03) (15, 3.6664e-03) (16, 3.3217e-03) (17, 3.0240e-03) (18, 2.7647e-03) (19, 2.6160e-03) (20, 2.3950e-03) (21, 2.2185e-03) (22, 2.0781e-03) (23, 1.9819e-03) (24, 1.8525e-03) (25, 1.7227e-03) (26, 1.6206e-03) (27, 1.5545e-03) (28, 1.4728e-03) (29, 1.4005e-03) (30, 1.3231e-03) (30, 1.3226e-03) (31, 1.2751e-03) (32, 1.2136e-03) (33, 1.1605e-03) (34, 1.1121e-03) (35, 1.0700e-03) (36, 1.0191e-03) (37, 9.8018e-04) (38, 9.4526e-04) (39, 9.1704e-04) (40, 8.7492e-04) 
			};
\addplot[mark=triangle, orange, line width = 0.5pt] coordinates {
				(1, 2.0232e-01)  (2, 1.4408e-01)  (3, 9.2207e-02)  (4, 7.5522e-02)  (5, 5.9411e-02)  (6, 5.1950e-02)  (7, 4.6086e-02)  (8, 4.0727e-02)  (9, 3.5747e-02)  (10, 3.2489e-02)  (11, 3.0145e-02)  (12, 2.7547e-02)  (13, 2.5566e-02)  (14, 2.4006e-02)  (15, 2.2286e-02)  (16, 2.0913e-02)  (17, 2.0084e-02)  (18, 1.8747e-02)  (19, 1.8048e-02)  (20, 1.7049e-02)  (21, 1.6205e-02)  (22, 1.5637e-02)  (23, 1.4977e-02)  (24, 1.4238e-02)  (25, 1.3761e-02)  (26, 1.3375e-02)  (27, 1.2808e-02)  (28, 1.2443e-02)  (29, 1.2076e-02)  (30, 1.1615e-02)  (31, 1.1235e-02)  (32, 1.0920e-02)  (33, 1.0566e-02)  (34, 1.0319e-02)  (35, 1.0087e-02)  (36, 9.6662e-03)  (37, 9.4791e-03)  (38, 9.2939e-03)  (39, 9.0644e-03)  (40, 8.8227e-03)  
			};
\addplot[mark=triangle*, orange, line width = 0.5pt] coordinates {
				(1, 8.2533e-01)  (2, 5.1041e-01)  (3, 4.5993e-01)  (4, 4.2320e-01)  (5, 3.0484e-01)  (6, 2.1470e-01)  (7, 2.3279e-01)  (8, 1.5763e-01)  (9, 1.4132e-01)  (10, 1.3966e-01)  (11, 1.4237e-01)  (12, 9.7238e-02)  (13, 9.8576e-02)  (14, 1.5781e-01)  (15, 9.4331e-02)  (16, 7.4903e-02)  (17, 7.7608e-02)  (18, 6.2756e-02)  (19, 6.2711e-02)  (20, 5.9514e-02)  (21, 5.9400e-02)  (22, 1.0343e-01)  (23, 5.9154e-02)  (24, 4.3502e-02)  (25, 4.0717e-02)  (26, 4.1174e-02)  (27, 4.0649e-02)  (28, 3.9182e-02)  (29, 3.8806e-02)  (30, 3.8978e-02)  (31, 3.7830e-02)  (32, 3.1358e-02)  (33, 3.3924e-02)  (34, 3.1513e-02)  (35, 3.1478e-02)  (36, 2.4346e-02)  (37, 2.3929e-02)  (38, 2.4475e-02)  (39, 2.5588e-02)  (40, 2.2727e-02)       
			};
\addplot[mark=diamond, line width = 0.5pt] coordinates {
				(1, 2.6463e-01)  (2, 1.1416e-01)  (3, 6.3187e-02)  (4, 4.4032e-02)  (5, 3.5890e-02)  (6, 2.4815e-02)  (7, 2.1730e-02)  (8, 1.6347e-02)  (9, 1.4136e-02)  (10, 1.1661e-02)  (11, 1.0911e-02)  (12, 8.4008e-03)  (13, 7.9015e-03)  (14, 7.0306e-03)  (15, 6.4736e-03)  (16, 5.4426e-03)  (17, 5.2489e-03)  (18, 4.5564e-03)  (19, 4.4058e-03)  (20, 3.8519e-03)  (21, 3.6036e-03)  (22, 3.4127e-03)  (23, 3.3241e-03)  (24, 2.7987e-03)  (25, 2.6895e-03)  (26, 2.5567e-03)  (27, 2.4528e-03)  (28, 2.2561e-03)  (29, 2.2094e-03)  (30, 1.9855e-03)  (31, 1.9581e-03)  (32, 1.7801e-03)  (33, 1.7258e-03)  (34, 1.6671e-03)  (35, 1.6208e-03)  (36, 1.4441e-03)  (37, 1.4207e-03)  (38, 1.3768e-03)  (39, 1.3448e-03)  (40, 1.2319e-03) 
			};
\addplot[mark=square, line width = 0.5pt] coordinates {
				(1, 4.4672e-01)  (2, 1.3226e-01)  (3, 1.1765e-01)  (4, 3.2729e-02)  (5, 3.5010e-02)  (6, 2.7088e-02)  (7, 2.7361e-02)  (8, 8.1672e-03)  (9, 8.0903e-03)  (10, 7.7766e-03)  (11, 7.8647e-03)  (12, 5.4141e-03)  (13, 5.4142e-03)  (14, 5.4144e-03)  (15, 5.3633e-03)  (16, 1.9248e-03)  (17, 1.9242e-03)  (18, 1.8076e-03)  (19, 1.8062e-03)  (20, 1.6773e-03)  (21, 1.6732e-03)  (22, 1.6585e-03)  (23, 1.6566e-03)  (24, 1.0192e-03)  (25, 1.0133e-03)  (26, 1.0055e-03)  (27, 1.0075e-03)  (28, 9.8244e-04)  (29, 9.8007e-04)  (30, 9.7061e-04)  (31, 9.7546e-04)  (32, 4.3734e-04)  (33, 4.3723e-04)  (34, 4.3456e-04)  (35, 4.3450e-04)  (36, 3.8048e-04)  (37, 3.7931e-04)  (38, 3.7950e-04)  (39, 3.7846e-04)  (40, 3.3234e-04)  
			};
\addplot[mark=diamond*, line width = 0.5pt] coordinates {
				(1, 8.2572e-01)  (2, 6.0007e-01)  (3, 5.4715e-01)  (4, 3.6068e-01)  (5, 3.4712e-01)  (6, 2.8834e-01)  (7, 2.8827e-01)  (8, 1.9877e-01)  (9, 1.8659e-01)  (10, 1.9106e-01)  (11, 1.9111e-01)  (12, 1.3127e-01)  (13, 1.3238e-01)  (14, 1.3276e-01)  (15, 1.3208e-01)  (16, 9.8024e-02)  (17, 9.8646e-02)  (18, 8.4433e-02)  (19, 8.4620e-02)  (20, 8.3351e-02)  (21, 8.1833e-02)  (22, 8.2821e-02)  (23, 8.2770e-02)  (24, 5.4217e-02)  (25, 5.5021e-02)  (26, 5.4195e-02)  (27, 5.3448e-02)  (28, 5.2997e-02)  (29, 5.3263e-02)  (30, 5.1943e-02)  (31, 5.2288e-02)  (32, 4.2222e-02)  (33, 4.2463e-02)  (34, 4.6917e-02)  (35, 4.2312e-02)  (36, 3.2904e-02)  (37, 3.2624e-02)  (38, 3.3223e-02)  (39, 3.2786e-02)  (40, 3.1663e-02) 
			};
\end{axis}
		\end{tikzpicture}
	\end{minipage}
	\caption{
		Comparing the approximation errors ${\varepsilon_2^R(h)}$ of the tensored $B_2$-cutoff \eqref{def:test_functions} approximated by various orthonormal systems in dimensions $d\in\{1,2,4,7\}$.
	}
	\label{fig:numeric_B2_ttc}
\end{figure}

\subsection{A note on $\ell_{\infty}$-approximation}
As derived in \cite{NaPo19} and recalled in \eqref{eq:periodization_idea}, the transformed Fourier system~\eqref{def:trafo_Fou_system} for non-periodic funtions is the result of applying an inverted change of variable $\psi^{-1}(\cdot,\bm\eta)$ in the form of~\eqref{def:Trafo_comb} to the Fourier system elements within the $L_{2}(\mathbb{T}^d)$-scalar product, in order to generate another orthonormal system in a given space $L_{2}\left([-\frac{1}{2},\frac{1}{2}]^d,\omega\right)$.
There are two interpretations for the resulting integral of the form
\begin{align}\label{eq:noname}
	\int_{\left[0,1\right]^d} \frac{\varrho(\mathbf y,\bm\eta)}{\omega(\mathbf y)}\,\mathrm{e}^{2\pi \mathrm{i}(\mathbf k-\mathbf m)\psi^{-1}(\mathbf y,\bm\eta)} \, \omega(\mathbf y) \,\mathrm{d}\mathbf y
	= \int_{\mathbb{T}^d} \mathrm{e}^{2\pi \mathrm{i}(\mathbf k-\mathbf m)\mathbf x} \,\mathrm{d}\mathbf x
	= \delta_{\mathbf k,\mathbf m}.
\end{align}
We either have another periodic system of the form $\left\{ \mathrm{e}^{2\pi\mathrm{i}\mathbf k\cdot\psi^{-1}(\cdot,\bm\eta)} \right\}_{\mathbf k\in I}$ and the weighted $L_{2}\left([0,1]^d,\varrho(\cdot,\bm\eta)\right)$-scalar product; 
or we attach $\sqrt{\varrho(\cdot,\bm\eta)/\omega(\cdot)}$ to the individual exponentials $\mathrm{e}^{2\pi\mathrm{i}\mathbf k\cdot\psi^{-1}(\cdot,\bm\eta)}$ and end up with the non-periodic system~\eqref{def:trafo_Fou_system} and the originally given weighted $L_{2}\left([-\frac{1}{2},\frac{1}{2}]^d,\omega\right)$-scalar product.
If we consider a constant weight function $\omega\equiv 1$, then there is a drawback that comes with the later choice, because $\varrho$ is unbounded and causes singularities at the boundary points of the elements in the approximated transformed Fourier sum~\eqref{def:approx_trafo_Fou_Partsum}.
So, the pointwise approximation error $\varepsilon_\infty^R$ in \eqref{def:ell_p_approxError} isn't finite, unless we consider a suitably weighted $\ell_\infty$-norm that counteracts the behavior of the approximant towards the boundary points, which is discussed more thoroughly in \cite{NaPo19}.
This strategy is based on choosing the weight function $\omega$ in such a way that the quotient $\varrho(\cdot,\bm\eta)/\omega(\cdot)$ is either constant or converges at the boundary points.
However, for any chosen torus-to-cube transformation - especially for the presented parameterized transformations $\psi(\cdot,\bm\eta)$ in \eqref{eq:logarithmic_trafo} and \eqref{eq:erf_trafo} with a fixed parameter $\bm\eta$ - the weight function has to be chosen in such a way so that on one hand the singularities of the density function are controlled and on the other hand the given function $h$ is still in $L_{2}(\left[0,1\right]^d,\omega)$.

We achieve this effect for example by showing the connection of the transformed Fourier framework with the Chebyshev system, when we put the Chebyshev transformation \eqref{def:Cheb_trafo} into the transformed Fourier system~\eqref{def:trafo_Fou_system} despite the fact that it is not a torus-to-cube transformation as in \eqref{def:Trafo_comb}.
Considering the hyperbolic cross $I_N^1$ as defined in \eqref{def:HC} and $x,y\in [0,1]$, we choose $\psi$ to be the Chebyshev transformation~\eqref{def:Cheb_trafo} of the form $\psi(x) = \frac{1}{2}+\frac{1}{2}\cos\left( 2\pi (x-\frac{1}{2}) \right)$, with the inverse $\psi^{-1}(y) = \frac{1}{2}+\frac{\arccos(2y-1)}{2\pi}$ and the density $\varrho(y) = \frac{1}{2\pi \sqrt{y(1-y)}}$.
By putting $\omega(y) = \varrho(y)$, the transformed Fourier system~\eqref{def:trafo_Fou_system} turns into
\begin{align}\label{eq:TrafoFourier_copy}
	\varphi_{k}(y)
	= \mathrm{e}^{\pi\mathrm{i}k + \mathrm i k\arccos(2y)}
	= (-1)^k(\cos(k\arccos(2y-1)) + \mathrm{i} \sin(k\arccos(2y-1)))
\end{align}
for $k\in\{-N,\ldots,N\}$ and by combining the positive and negative frequencies we obtain
\begin{align*}
	\varphi_{k}(y) = 
	\begin{cases}
		1 \quad &\text{ for } \quad k = 0,\\
	 	(-1)^k 2\cos(k\arccos(2y-1)) \quad &\text{ for } \quad k \in\{1,2,\ldots,N\},
	\end{cases}
\end{align*}
which is orthogonal with respect to the $L_2\left([0,1], \omega \right)$-scalar product with $\omega(y) = \frac{1}{2\pi \sqrt{y(1-y)}}$.
With some additional scaling we obtain an orthonormal system that's equivalent to the Chebyshev system~\eqref{def:T_k}.

\section{Conclusion}
We considered the approximation of non-periodic functions on the cube $\left[0,1\right]^d$ by different systems of orthonormal functions.
We compared the Chebyshev system that is orthonormal with respect to a weighted $L_2$-scalar product, 
the system of half-periodic cosines that uses tent-transformed sampling nodes
and a parameterized transformed Fourier system.
For the cosine system, which basically only mirrors a non-periodic function at it's boundary points, as well as the transformed Fourier system with a small parameter, yielded the worst approximation errors.
Switching to the Chebyshev system, which mirrors and additionally smoothens a given function, improved the approximation error decay.
The same effect was obtained for the transformed Fourier system after increasing the parameter enough to obtain a better smoothening effect.
The numerical experiments showcased the proposed parameter control in \cite{NaPo19} that is set up by periodizing functions via families of parameterized torus-to-cube mappings.
This approach in particular generalizes the idea used to derive Chebyshev polynomials.

\section*{Acknowledgements}
The authors thank the referees for their valuable suggestions and remarks.
The first named author gratefully acknowledges the support by the funding of the European Union and the Free State of Saxony (ESF).

\small
\bibliographystyle{abbrv}

\end{document}